\pdfoutput=1
\documentclass[journal]{IEEEtran}
\usepackage{amsfonts}
\usepackage{bm}
\usepackage{amsmath,amsfonts,amsthm,amssymb}
\usepackage{setspace}
\usepackage{Tabbing}
\usepackage{fancyhdr}
\usepackage{lastpage}
\usepackage{extramarks}
\usepackage{chngpage}
\usepackage{algorithmic}
\usepackage{soul,color}
\usepackage{epsfig}
\usepackage{graphicx,float,wrapfig,subfigure}
\usepackage{dsfont}
\usepackage{longtable}
\usepackage{wrapfig}
\usepackage{stfloats}
\usepackage{cite}
\usepackage{graphicx}
\usepackage{array}
\usepackage{multirow}{\tiny }
\usepackage{multicol}
\usepackage{colortbl}
\usepackage{tabularx}
\usepackage{mdwmath}
\usepackage{mdwtab}
\usepackage{color}
\usepackage{verbatim}
\usepackage{amsmath}
\usepackage{tikz}
\usetikzlibrary{calc}
\usepackage{flowchart}
\usetikzlibrary{shapes.geometric, arrows}
\usepackage{overpic}
\usepackage{epstopdf}
\usepackage{url}
\usepackage[colorlinks,linkcolor=black,anchorcolor=black,citecolor=black,urlcolor=black]{hyperref} 

\usepackage{makecell}
\usepackage{textcomp}
\usepackage{algorithm}
\usepackage{algorithmic}
\usepackage{threeparttable}

\newtheorem{proposition}{Proposition}
\newtheorem{lemma}{Lemma}

\newtheorem{remark}{Remark}

\usepackage{color}
\makeatletter %
\let\myorg@bibitem\bibitem
\def\bibitem#1#2\par{%
	\@ifundefined{bibitem@#1}{%
		\myorg@bibitem{#1}#2\par
	}{%
		\begingroup
		\color{\csname bibitem@#1\endcsname}%
		\myorg@bibitem{#1}#2\par
		\endgroup
	}%
}

\makeatother %

\allowdisplaybreaks

\usepackage{footnote}
\makesavenoteenv{tabular}
\makesavenoteenv{table}

\begin{document}

\title{
Chance-constrained regulation capacity offering for HVAC systems under non-Gaussian uncertainties with mixture-model-based convexification
	}

\author{
Ge~Chen,~\IEEEmembership{Graduate Student Member,~IEEE,}
Hongcai~Zhang,~\IEEEmembership{Member,~IEEE,}
Hongxun~Hui,~\IEEEmembership{Member,~IEEE,}
and~Yonghua~Song,~\IEEEmembership{Fellow,~IEEE}

\thanks{	
	G. Chen, H. Zhang, H. Hui, and Y. Song are with the State Key Laboratory of Internet of Things for Smart City and Department of Electrical and Computer Engineering, University of Macau, Macao, 999078 China (email: hczhang@um.edu.mo).
	}
}

\maketitle

\begin{abstract}
Heating, ventilation, and air-conditioning (HVAC) systems are ideal demand-side flexible resources to provide regulation services. However, finding the best hourly regulation capacity offers for HVAC systems in a power market ahead of time is challenging because they are affected by non-Gaussian uncertainties from regulation signals. Moreover, since HVAC systems need to frequently regulate their power according to regulation signals, numerous thermodynamic constraints are introduced, leading to a huge computational burden. This paper proposes a tractable chance-constrained model to address these challenges. It first develops a temporal compression approach, in which the extreme indoor temperatures in the operating hour are estimated and restricted in the comfortable range so that the numerous thermodynamic constraints can be compressed into only a few ones. Then, a novel convexification method is proposed to handle the non-Gaussian uncertainties. This method leverages the Gaussian mixture model to reformulate the chance constraints with non-Gaussian uncertainties on the left-hand side into deterministic non-convex forms. We further prove that these non-convex forms can be approximately convexified by second-order cone constraints with 
marginal optimality loss. Therefore, the proposed model can be efficiently solved with guaranteed optimality. Numerical experiments are conducted to validate the superiority of the proposed method.

\end{abstract}

\begin{IEEEkeywords}
HVAC systems, demand-side flexibility, regulation capacity, chance-constrained programming, Gaussian mixture model, convexification.
\end{IEEEkeywords}

\section{Introduction}
The growing penetration of renewable energies in power systems reduces fossil fuel consumption and carbon emissions. However, the intermittent and stochastic characteristics of renewable energies may cause power fluctuation problems, which severely threatens the stability of power systems \cite{7738432}. To support the stable operation of power systems, more flexible resources are required to participate in regulation services \cite{impram2020challenges}.

Heating, ventilation, and air conditioning (HVAC) loads are one of the most promising demand-side flexible resources because of the building's inherent ability to store heating/cooling power \cite{kohlhepp2019large,8295258}. To utilize the flexible HVAC systems for regulation services, the corresponding regulation capacity offers need to be reported to the power market in advance \cite{dispatch2019pjm}. The power market can collect all the regulation capacity offers to design regulation signals. Then, HVAC systems can follow regulation signals and adjust their power scheduling to earn regulation revenue \cite{7527683, 8648288}. Since this revenue is in proportion to the regulation capacity offers, increasing attention has been paid to quantifying the potential regulation capacity for HVAC systems. For example, reference \cite{7864461} proposed a geometric approach to characterize the aggregated regulation capacity of HVAC systems. Reference \cite{8731716} proposed a robust-based method to quantify the HVAC's regulation capacity in distribution networks. Reference \cite{9502573} leveraged deep learning techniques to develop a model-free method to determine the best regulation capacity for HVAC systems. Because one HVAC system's thermal inertia is limited, its regulation capacity can get affected by uncertain and biased regulation operations. However, most of the aforementioned papers do not consider the impacts of regulation signals, which may overly estimate HVAC's regulation capacity and violate the corresponding building's indoor thermal comforts.

Nevertheless, taking regulation signals into consideration is challenging because they are highly stochastic and unpredictable. To address this issue, some papers treat signals as uncertainties and leverage robust optimization to design perfectly safe scheduling strategies for flexible resources \cite{7422803,8168425}. However, robust optimization methods do not allow any constraint violation for all realizations of uncertainties \cite{bertsimas2018data}, so their solutions are usually overly conservative. An alternative choice is chance-constrained programming (CCP). CCP only requires constraints to be satisfied with a predetermined probability and allow small violations so that it can better balance robustness and optimality \cite{geng2019data,9535415}. Considering that little thermal discomfort can be temporarily tolerated, CCP is preferable for optimizing the HVAC's regulation capacity.
Unfortunately, applying CCP still faces two challenges:
{\begin{enumerate}
\item {Regulation signals do not follow Gaussian distribution. However, 
the most commonly used CCP method is based on the Gaussian assumption, in which the uncertainty is assumed to be normally distributed \cite{8017474,9122389,9535415}. Thus, applying this Gaussian-assumption-based method may lead to infeasible solutions. Some other scholars proposed distributionally robust chance-constrained methods (DRCC) to handle non-Gaussian uncertainties. Based on specific moment information or Wasserstein distance, DRCC constructs an ambiguity set to cover possible distributions and requires the probabilistic constraints to be robust to the ambiguity set \cite{esfahani2018data,9417102}. DRCC can handle non-Gaussian uncertainties and has been used for scheduling flexible sources under uncertain regulation signals \cite{8301597,9347818}.
However, DRCC may still result in overly conservative solutions and dramatically reduce the regulation revenue because the ambiguity set may cover some distributions that are much different from the true one.}
\item The regulation signal updates frequently (e.g. RegD signal updates every two seconds in the PJM market \cite{dispatch2019pjm}). Since the HVAC power needs to be regulated to follow these signals, numerous constraints are involved in order to ensure thermal comforts corresponding to each regulation signal. This will make the CCP computationally expensive, especially for the methods that involve many additional variables and constraints for each chance constraint (e.g. sample average approximation \cite{8680681}, conditional value-at-risk approximation \cite{7828060}).
\end{enumerate}}

To overcome the first challenge, several papers combined the Gaussian mixture model (GMM) with CCP. GMM is a universal approximator of probability densities, and any non-Gaussian distribution can be approximately fitted with a finite number of Gaussian components \cite{goodfellow2016deep}. 
{In \cite{8772186,8936474}, GMM was used to fit the non-Gaussian renewable energy uncertainties. Then, the chance constraints were directly reformulated into tractable forms based on the quantile of uncertainties.}
In \cite{9376652}, an online-offline double-track approach was developed to accelerate the fitting of GMM for the uncertainties of gas demands. {However, these GMM-based methods are only suitable for the chance constraints with right-hand side (RHS) uncertainties\footnote{{Consider a linear constraint $\bm a^\intercal \bm x \leq b$. If the uncertainty is the vector $\bm a$, then we call it ``left-hand side (LHS) uncertainty"; if the uncertainty is in the constant $b$, then it is called ``right-hand side (RHS) uncertainty".}}. Considering that the regulation signal uncertainties are on the left-hand side (LHS), these methods are still inapplicable.}


Unlike the first challenge, there are only very few papers that have tried to tackle the second challenge. In fact, most published papers, including \cite{7422803,8168425,9347818}, only required the satisfaction of constraints with low temporal resolutions, while the intermediate variable variations between two neighboring time slots were ignored. Thus, although this manner can reduce the constraint number, it cannot always guarantee feasibility. For example, reference \cite{9347818} only restricted that constraints should be satisfied every five minutes, e.g., at $t \in$\{5min, 10min, $\cdots$, 55min, 60min\}. However, it cannot guarantee that there is no violation within each five minutes, e.g., at $t$=7min.

To overcome the aforementioned two challenges, we propose a tractable chance-constrained model to optimize the regulation capacity offering for HVAC systems. The specific contributions are threefold:
\begin{enumerate}
    \item We propose a chance-constrained model to determine the hour-ahead regulation capacity offers for HVAC systems in the PJM market. This model considers the impacts of non-Gaussian uncertainties from regulation signals. Moreover, the thermodynamic constraints are built according to the updated frequency of regulation signals (i.e. every two seconds) so that indoor thermal comforts can be properly maintained. 
	\item To address the intractability from non-Gaussian uncertainties on the LHS, we propose a mixture-model-based convexification method. It first leverages GMM to reformulate each chance constraint with non-Gaussian LHS uncertainties into a deterministic non-convex form. Then, this non-convex reformulation is equivalently re-expressed as an exponential form. Based on piece-wise linearization, we further prove that this exponential form can be safely approximated by an SOCP constraint with marginal optimality loss, which guarantees desirable optimality and computational efficiency. 
	To the best of our knowledge, this is the first time that GMM-based methods can be extended to chance constraints with LHS uncertainties from regulation signals.
    \item To reduce the computational burden brought by the huge number of thermodynamic constraints, we propose a temporal compression method. In this method, we first estimate the maximum and minimum indoor temperatures over a long time duration based on the monotonicity of the thermodynamic model. Then, by restricting the estimated extreme temperatures in the comfortable range, all thermodynamic constraints in this long time duration can be replaced by only a few ones, which significantly enhances computational efficiency.
\end{enumerate}


The remaining parts are organized as follows. Section \ref{sec_formula} describes the problem formulation. Section \ref{sec_solution} presents the details of the proposed mixture-model-based convexification method. Section \ref{sec_case} shows simulation results and Section \ref{sec_conclusion} concludes this paper.

\section{Problem formulation} \label{sec_formula}

\begin{figure}
	\subfigbottomskip=-6pt
	\subfigcapskip=-4pt
		\vspace{-4mm}
	\centering
	{\includegraphics[width=0.9\columnwidth]{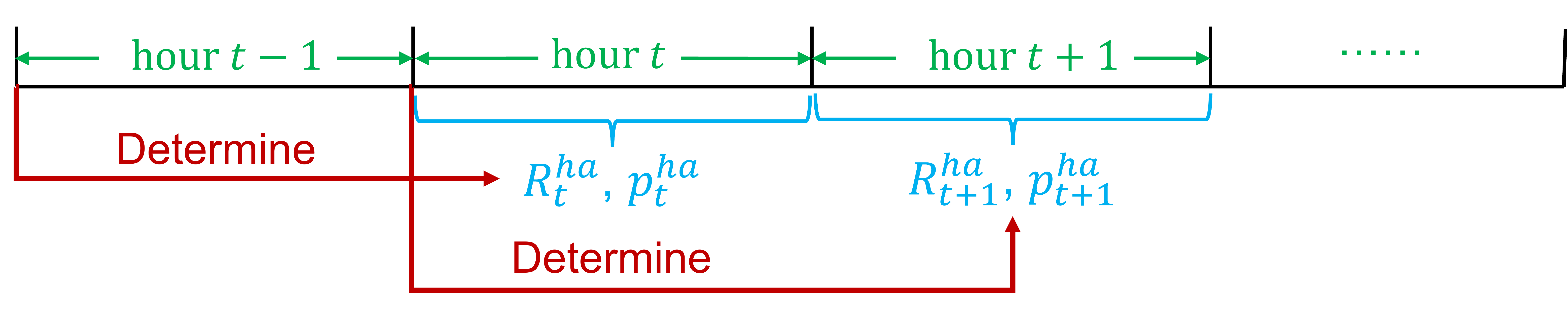}}
	\vspace{-4mm}
	\caption{Schema of hour-ahead regulation capacity offers in the PJM market.}
	\label{fig_market}
	\vspace{-4mm}
\end{figure}

{We consider an aggregator strategically operating a couple of HVAC systems to provide regulation capacities in the PJM market. As required, the aggregator shall offer the regulation capacity to the market at least one hour ahead, as shown in Fig. \ref{fig_market}. For example, the capacity offer for 3:00pm-4:00pm should be reported before 2:00pm. In order to maximize its regulation revenue, the aggregator need to properly design the power schedule of HVAC systems, i.e., $p_t^\text{ha}$, and accurately estimate their corresponding regulation capacities, i.e., $R_t^\text{ha}$.
Because buildings have limited thermal inertia, their regulation capacities are significantly affected by the uncertain regulation signals, which should be explicitly considered.
}

Since the regulation signal in PJM updates every two seconds, we divide the operating hour into 1800 time slots with $l \in \mathcal{L}$ as their indexes and $\Delta L=2$ seconds as the time step size. By using $i \in \mathcal{I}$ to index HVAC systems, the thermodynamic model of one building can be expressed as:
\begin{align}
C_i \frac{d\theta_i^\text{in}}{dt} = g(\theta^\text{out}-\theta_i^\text{in}) + h_i - q_i, \forall i \in \mathcal{I}, \label{eqn_thermal_diff}
\end{align}
where $C_i$ and $g_i$ are the building heat capacity and heat transfer coefficient between indoor and outdoor environments, respectively; $\theta^\text{in}_i$ and $\theta^\text{out}$ are the indoor and outdoor temperatures, respectively; $h_i$ and $q_i$ denote the heat load contributed by indoor sources (e.g. electronic devices) and cooling supply from HVAC systems, respectively. 

{Since the outdoor temperature $\theta^\text{out}$ and indoor heat load $h_i$ vary slowly, we assume that they keep unchanged in the operating hour.} When the HVAC system participates in regulation services, their power need to respond to the regulation signals. Thus, the actual HVAC power at the $l$-th time slot in hour $t$, i.e., $p_{t, l,i}^\text{HV}$, and corresponding cooling supply, i.e., $q_{t, l,i}$, are expressed as: 
\begin{align} 
	&p_{t,l,i}^\text{HV} = p_{t,i}^\text{ha} - R_{t,i}^\text{ha} s_{t,l},\quad \forall  l \in \mathcal{L}, \forall i \in \mathcal{I}, \label{eqn_p}\\
	&q_{t,l,i}=\text{COP}_i \cdot p_{t,l,i}^\text{HV}, \quad \forall  l \in \mathcal{L}, \forall i \in \mathcal{I}, \label{eqn_q}
\end{align}
where $s_{t,l}$ is the $l$-th signal in hour $t$; $\text{COP}_i$ is the coefficient of performance of the HVAC system in building $i$. 
Eqs. (\ref{eqn_p})-(\ref{eqn_q}) indicate that the cooling supply keeps constant in each $\Delta l$, {so  (\ref{eqn_thermal_diff}) can be directly integrated from time $l$ to $l+\Delta l$: }
\begin{align}
\theta_{t,l,i}^\text{in}=&a_i^\text{in}\theta_{t,l-1,i}^\text{in}+a_i^\text{out}\theta_{t}^\text{out} + a_i^\text{h}h_t \notag\\
&+ a_i^\text{q}(p_{i,t}^\text{ha} - R_{t,i} s_{t,l-1}), \quad \forall l \in \mathcal{L}/\{0\}, \forall i \in \mathcal{I}, \label{eqn_thermal}
\end{align}
where $a_i^\text{in}$, $a_i^\text{out}$, $a_i^\text{h}$, and $a_i^\text{q}$ are thermal parameters, which are calculated based on $C_i$, $g_i$, and $\text{COP}_i$ in  (\ref{eqn_thermal_diff})-(\ref{eqn_q}). Since $s_{t, l}$ is uncertain, both the actual HVAC power $p_{t, l,i}^\text{HV}$ and indoor temperature $\theta_{t,l,i}^\text{in}$ are also uncertain according to  (\ref{eqn_p})-(\ref{eqn_thermal}).  Considering temporary thermal discomfort in extreme conditions can be tolerated, we employ chance-constrained programming to describe the thermal comfort requirement:
\begin{align}
\begin{cases}
&\mathbb{P}\left( \theta_{t,l,i}^\text{in} \leq \theta_i^\text{max}\right)\geq 1-\epsilon,\\
&\mathbb{P}\left( \theta_{t,l,i}^\text{in} \geq \theta_i^\text{min} \right)\geq 1-\epsilon, 
\end{cases} \quad \forall l \in \mathcal{L}, \quad \forall i \in \mathcal{I}, \label{eqn_cc_o1}
\end{align}
where $\theta^\text{max}$ and $\theta^\text{min}$ are the upper and lower bounds of the thermal comfortable range, respectively; $\epsilon$ is the risk parameter. The device limit requires that the HVAC power should always stay in the allowable range, so a robust manner is used to describe this limitation:
\begin{align}
\begin{cases}
&\max_{s_{t, l}} p_{t,l,i}^\text{HV} \leq p_i^\text{max}, \\
&\min_{s_{t, l}} p_{t,l,i}^\text{HV} \geq p_i^\text{min},
\end{cases} \quad \forall l \in \mathcal{L}, \quad \forall i \in \mathcal{I}, \label{eqn_cc_o2}
\end{align}
where $p_i^\text{max}$ and $p_i^\text{min}$ are the upper and lower bounds of the HVAC power, respectively. Since the regulation signal is restricted in [-1,1], Eq. (\ref{eqn_cc_o2}) can be reformulated as follows according to  (\ref{eqn_p}):
\begin{align}
p_{t,i}^\text{ha} + R_{t,i}^\text{ha} \leq p_i^\text{max}, \quad p_{t,i}^\text{ha} - R_{t,i}^\text{ha} \geq p_i^\text{min}, \quad \forall i \in \mathcal{I}. \label{eqn_r_power}
\end{align}
In the PJM market, the hour-ahead regulation capacity offers are only allowed to be reduced from the the corresponding day-ahead offers, i.e., $R_{t,i}^\text{da}$, as follows:
\begin{align}
0 \leq R_{t,i}^\text{ha} \leq R_{t,i}^\text{da}, \quad \forall i \in \mathcal{I}, \label{eqn_R_bound}
\end{align}
where $R_{t,i}^\text{da}$ is determined one day ahead, so it is a known parameter when we optimize the hour-ahead offer $R_{t,i}^\text{ha}$.

Our objective is to minimize the total cost $EC_t$, which equals to the energy cost minus the regulation revenue:
\begin{align}
    EC_t =  \sum_{i} \sum_{l \in \mathcal{L}}\eta_t p_{t,i,l}^\text{HV}\Delta l -\sum_{i}(r_t^\text{rc}+r_t^{m}m_t)R_{t,i}^\text{ha}, \label{eqn_obj} 
\end{align}
{where the first and second terms on the RHS of  (\ref{eqn_obj}) represent the energy cost of HVAC systems and revenue from the HVAC power regulation.} Symbol $\mathbb{E}$ denotes the expectation operator; $\eta_t$ is the price for electricity purchasing at hour $t$; $\Delta t=1h$ is the operating time duration; $r_t^\text{rc}$ and $r_t^{m}$ are the unit revenues for regulation capacity and millage, respectively; $m_t$ is the regulation millage, which is defined as:
\begin{align}
    m_t = \sum_{l} | s_{t,l+1}-s_{t}|. \label{eqn_millage}
\end{align}
Considering $s_t$ is uncertain regulation signals, the regulation millage $m_t$ is also uncertain in this optimization problem.

Finally, the optimization problem is formulated as:
\begin{align} 
	&\min_{R_{t,i}^\text{ha}, p_{t,i}^\text{ha}, \forall l, \forall i} \mathbb E (EC_t) \tag{$\textbf{P1}$},\text{ s.t.:} \text{ Eqs.} \text{ (\ref{eqn_p})-(\ref{eqn_cc_o1}) and (\ref{eqn_r_power})-(\ref{eqn_obj}).}
\end{align}

{Solving \textbf{P1} is quite challenging. 
On the one hand, due to the high update frequency of regulation signals, Eq. (\ref{eqn_cc_o1}) introduces numerous thermodynamic constraints, leading to computational intractability. On the other hand, in  (\ref{eqn_cc_o1}), the signal uncertainties are on the LHS and do not follow the Gaussian distribution. Therefore, not only the Gaussian-assumption-based models used in \cite{8017474,9122389,9535415} but also the GMM-based methods proposed in \cite{8772186,8936474,9376652} can not be directly applied.}

\section{Solution Methodology} \label{sec_solution}
To overcome the aforementioned challenges, we first propose a temporal compression approach to reduce the thermodynamic constraint number. Then, a mixture-model-based convexification method is developed to reformulate the chance constraint with non-Gaussian uncertainties, i.e., Eq. (\ref{eqn_cc_o1}), into second-order cone programming (SOCP) forms. For simplicity, we omit the subscripts $t$ and $i$ in this section.

\subsection{Temporal compression} \label{sec_TC}
The key idea of the proposed temporal compression is to estimate the maximum and minimum indoor temperatures in the operating hour. Then, by restricting these extreme indoor temperatures in the comfortable range, the thermodynamic constraint number in  (\ref{eqn_cc_o1}) can be significantly reduced. For example, the first line in  (\ref{eqn_cc_o1}), i.e., $\mathbb{P}\left( \theta_{l}^\text{in} \leq \theta^\text{max}\right)\geq 1-\epsilon, \forall l \in \mathcal{L}$, contains 1800 constraints. However, once the maximum indoor temperature in the operating hour, i.e., $\overline \theta = \max_{l \in \mathcal{L}} \theta_{l}$ is estimated, the 1800 constraints can be replaced by only one single constraint $\mathbb{P}\left(\overline \theta \leq \theta^\text{max}\right)\geq 1-\epsilon$. Thus, the key problem is how to estimate these maximum and minimum indoor temperatures.

According to  (\ref{eqn_thermal}), the $l$-th indoor temperature in hour $t$, i.e., $\theta_{l}^\text{in}$, is expressed as:
\begin{align}
\theta_{l}^\text{in}=&a^\text{in}\theta_{l-1}^\text{in}+a^\text{out}\theta^\text{out} + a^\text{h}h + a^\text{q}(p^\text{ha} - R^\text{ha}s_{l-1}) \notag \\
=& a^\text{in}\left(a^\text{in}\theta_{l-2}^\text{in}+a^\text{out}\theta^\text{out} + a^\text{h}h+ a^\text{q} (p^\text{ha} - R^\text{ha}s_{l-2})\right) \notag\\ 
&+ a^\text{out}\theta^\text{out} + a^\text{h}h + a^\text{q}(p^\text{ha} - R^\text{ha}s_{l-1}) \notag\\
=& \cdots\cdots \notag \\
=& R^\text{ha}[\bm A \bm s]_{l} + f(l), \quad \forall l \in \mathcal{L}/\{0\},\label{eqn_thermal_tc}
\end{align}
where $\bm A$ is a coefficient matrix calculated based on $a^\text{q}$ and $a^\text{in}$;
vector $\bm s$ represents $[s_0, s_1, \cdots, s_{|\mathcal{L}|-1}]$; $[\bm A \bm s]_{l}$ denotes the $l$-th element of the product $\bm A \bm s$; function $f(l)$ is defined as:
\begin{align}
f(l) =& (a^\text{in})^{l} \theta_{0}^\text{in} + a^\text{out} \frac{1-(a^\text{in})^{l}}{1-a^\text{in}}\theta^\text{out} \notag \\
&+ a^\text{h} \frac{1-(a^\text{in})^{l}}{1-a^\text{in}}h + a^\text{q} \frac{1-(a^\text{in})^{l}}{1-a^\text{in}} p^\text{ha}. \label{eqn_f}
\end{align}

\begin{figure}
	\centering
	 			\vspace{-4mm}
	\includegraphics[width=1\columnwidth]{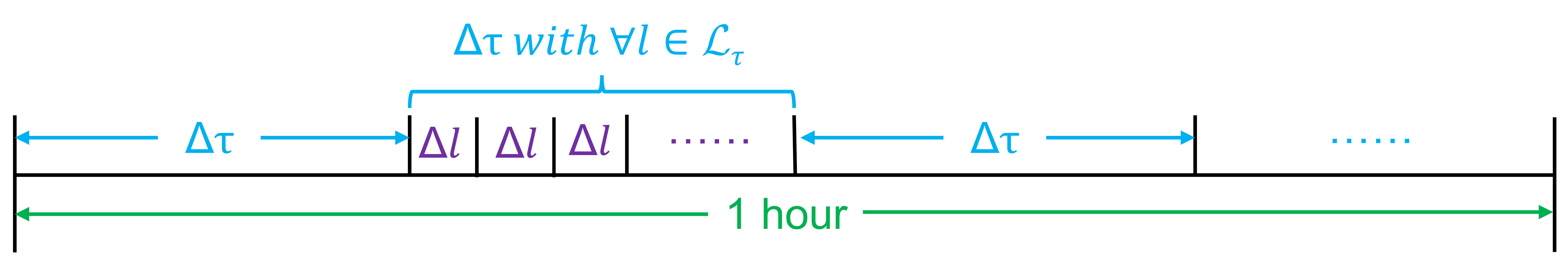}\vspace{-4mm}
	\caption{Schematic diagram for uniformly splitting the operating hour into multiple shorter duration $\Delta \tau$ , where $\Delta l=2s$ denotes the update interval of regulation signals; set $\mathcal{L}_{\tau}$ contains the indexes of $l$ in the $\tau$-th time duration.  
	}
	\label{fig_TC}
	 		\vspace{-4mm}
\end{figure}

Based on  (\ref{eqn_thermal_tc}), we must have:
\begin{align}
\begin{cases}
&\max_{l \in \mathcal{L}} \theta_{l}^\text{in} \leq \max_{l \in \mathcal{L}} f(l) + R \max_{l \in \mathcal{L}}[\bm A \bm s]_{l}, \\
&\min_{l \in \mathcal{L}} \theta_{l}^\text{in} \geq \min_{l \in \mathcal{L}} f(l) + R \min_{l \in \mathcal{L}}[\bm A \bm s]_{l}.
\end{cases} \label{eqn_approximation_0}
\end{align}
Thus, the maximum and minimum indoor temperatures can be approximated by the RHS terms of (\ref{eqn_approximation_0}). However, since $\max_{l \in \mathcal{L}} f(l)$ and $\max_{l \in \mathcal{L}}[\bm A \bm s]_{l}$ may appear at different moments, Eq. (\ref{eqn_approximation_0}) may lead to overly conservative solutions. To mitigate this conservativeness, we uniformly split the operating hour $\Delta t=1h$ into multiple shorter time duration $\Delta \tau$, as shown in Fig. \ref{fig_TC}. Then, the extreme indoor temperatures in each $\Delta \tau$ can be approximated by the RHS terms of (\ref{eqn_approximation}):
\begin{align}
\begin{cases}
&\max_{l \in \mathcal{L}_\tau} \theta_{l}^\text{in} \leq \max_{l \in \mathcal{L}_\tau} f(l) + R \max_{l \in \mathcal{L_\tau}}[\bm A \bm s]_{l}, \\
&\min_{l \in \mathcal{L}_\tau} \theta_{l}^\text{in} \geq \min_{l \in \mathcal{L}_\tau} f(l) + R \min_{l \in \mathcal{L}_\tau}[\bm A \bm s]_{l}, \notag
\end{cases}\\
\quad \quad\quad\quad\quad \forall \tau \in \mathcal{T}, \label{eqn_approximation}
\end{align}
where $\tau \in \mathcal{T}$ is the index of the shorter duration; $\mathcal{L}_{\tau}$ denotes the index set of $l$ in the $\tau$-th duration, which is obtained by uniformly splitting $\mathcal{L}$ into $|\mathcal{T}|$ parts. By adding a maximum/minimum operator over $\tau \in \mathcal{T}$ on both sides of (\ref{eqn_approximation}), we have: 
\begin{align}
\begin{cases}
\max_{l \in \mathcal{L}} \theta_{l}^\text{in} \leq \max_{\tau \in \mathcal{T}} \{\max_{l \in \mathcal{L}_\tau} f(l) + R \max_{l \in \mathcal{L_\tau}}[\bm A \bm s]_{l}\}, \\
\min_{l \in \mathcal{L}} \theta_{l}^\text{in} \geq \min_{\tau \in \mathcal{T}} \{\min_{l \in \mathcal{L}_\tau} f(l) + R \min_{l \in \mathcal{L_\tau}}[\bm A \bm s]_{l}\}.
\end{cases} \label{eqn_approximation_2}
\end{align}
Then, the extreme indoor temperatures in the operating hour can be estimated by the RHS terms in (\ref{eqn_approximation_2}).

\begin{proposition} \label{prop_0}
The approximation in (\ref{eqn_approximation_2}) is less conservative than that in (\ref{eqn_approximation_0}).
\end{proposition}
\emph{Proof}: See Appendix \ref{app_0}. 

Based on  (\ref{eqn_f}), function $f(l)$ is monotone with respect to $l$ because $\theta_{0}^\text{in}$, $\theta^\text{out}$, $h$, and $p^\text{ha}$ keep unchanged in the operating hour. Thus, the maximum and minimum values of $f(l)$ in each $\Delta \tau$ must appear at the boundaries, i.e., at $l=\tau\frac{|\mathcal{L}|}{|\mathcal{T}|}$ or $l= (\tau+1)\frac{|\mathcal{L}|}{|\mathcal{T}|}$. As for the terms $\max_{l \in \mathcal{L}_\tau}[\bm A \bm s]_{l}$ and $\min_{l \in \mathcal{L}_\tau}[\bm A \bm s]_{l}$, their values are uncertain but independent of decision variables. Therefore, we can directly treat them as two new uncertain parameters, i.e., $\overline u_\tau$ and $\underline u_\tau$, as follows:
\begin{align}
\overline u_\tau = \max_{l \in \mathcal{L}_\tau}[\bm A \bm s]_{l}, \quad \underline u_\tau = \min_{l \in \mathcal{L}_\tau}[\bm A \bm s]_{l}, \quad \forall \tau \in \mathcal{T}. \label{eqn_u_extreme}
\end{align}
Finally, Eq. (\ref{eqn_cc_o1}) can be replaced by:
\begin{align}
\begin{cases}
&\mathbb{P}\left(f\left(\tau\frac{|\mathcal{L}|}{|\mathcal{T}|}\right) + \overline u_\tau R^\text{ha} \leq \theta^\text{max} \right)\geq 1-\epsilon,\\
&\mathbb{P}\left(f\left((\tau+1)\frac{|\mathcal{L}|}{|\mathcal{T}|}\right) + \overline u_\tau R^\text{ha} \leq \theta^\text{max} \right)\geq 1-\epsilon,\\
&\mathbb{P}\left(f\left(\tau\frac{|\mathcal{L}|}{|\mathcal{T}|}\right) + \underline u_\tau R^\text{ha} \geq \theta^\text{min} \right)\geq 1-\epsilon,  \\
&\mathbb{P}\left(f\left((\tau+1)\frac{|\mathcal{L}|}{|\mathcal{T}|}\right) + \underline u_\tau R^\text{ha} \geq \theta^\text{min} \right)\geq 1-\epsilon, 
\end{cases}\quad \forall \tau \in \mathcal{T}. \label{eqn_thermal_r}
\end{align}
Based on (\ref{eqn_approximation_2}), any feasible solution of (\ref{eqn_thermal_r}) must be also feasible for the chance constraint (\ref{eqn_cc_o1}). Thus, Eq. (\ref{eqn_thermal_r}) is a safe approximation of (\ref{eqn_cc_o1}). Moreover, the thermodynamic constraint number is reduced from $2 \cdot |\mathcal{L}|$ to $4 \cdot |\mathcal{T}|$ (note that $|\mathcal{T}| \ll |\mathcal{L}|$), which significantly reduces the computational burden. 

\begin{figure}
	\subfigbottomskip=-6pt
	\subfigcapskip=-4pt
	\centering
	\subfigure[]{\includegraphics[width=0.49\columnwidth]{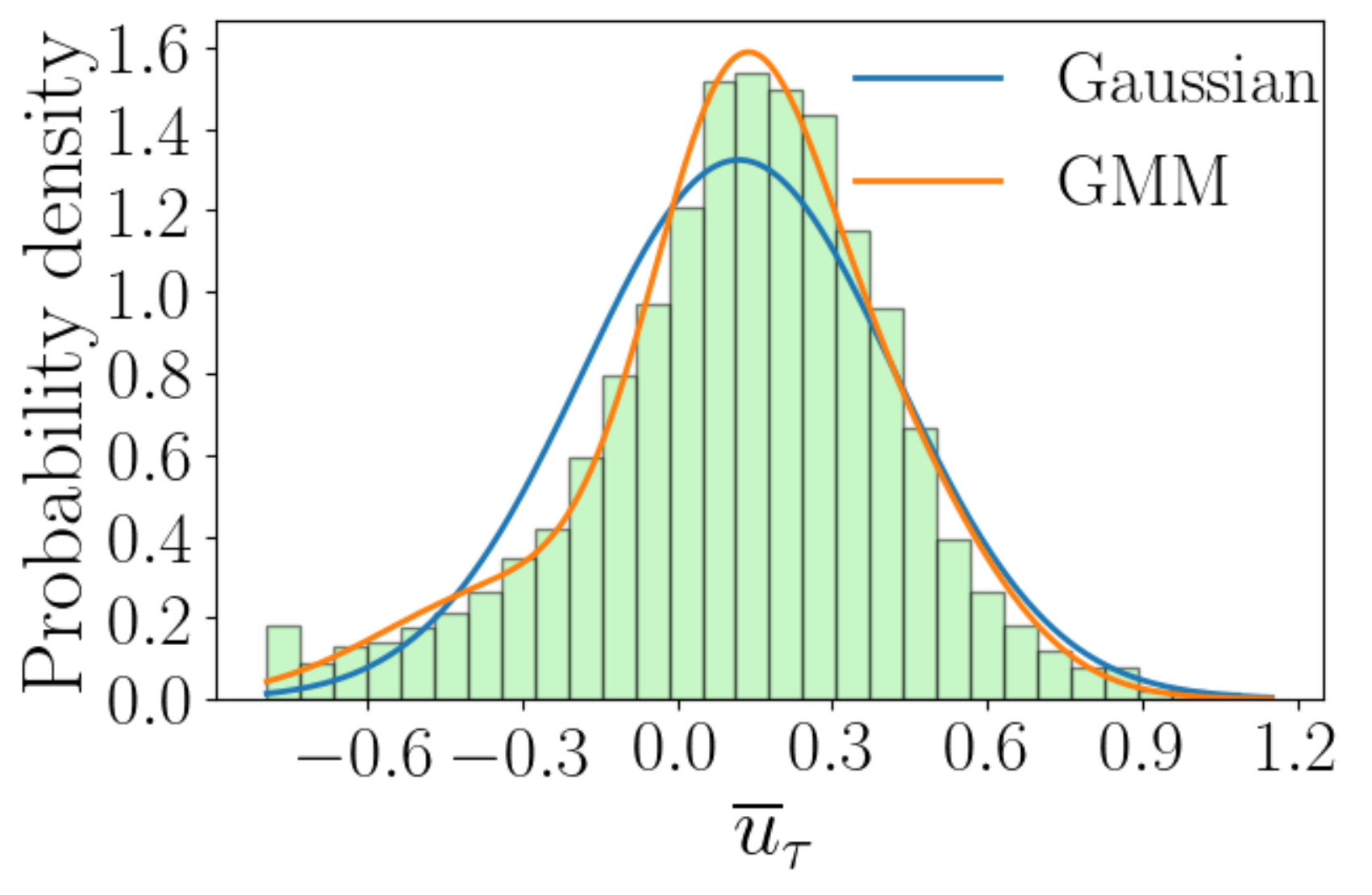}}
	\subfigure[]{\includegraphics[width=0.49\columnwidth]{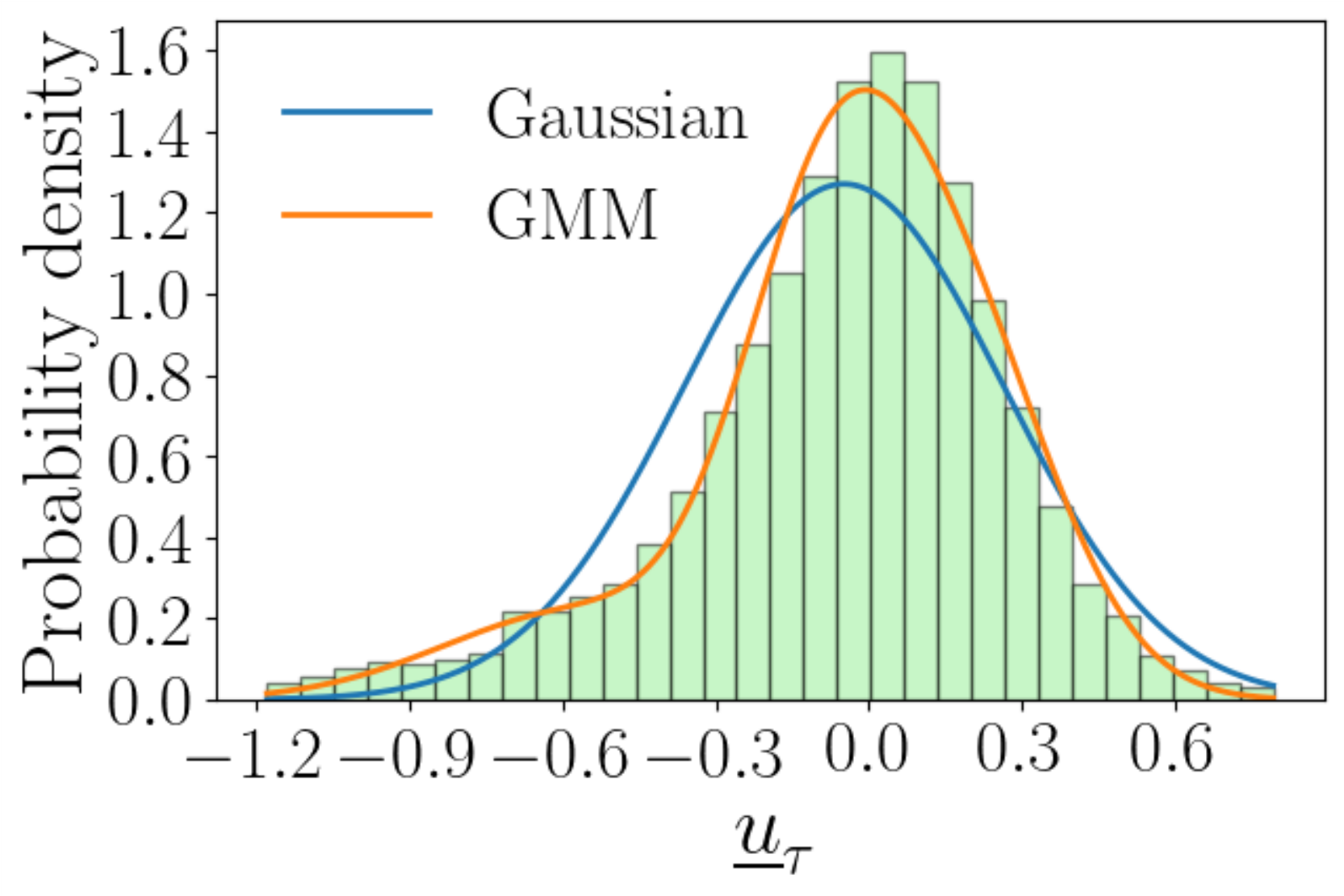}}
	\vspace{-4mm}
	\caption{The probability distribution of (a) $\overline u_\tau$ and (b) $\underline u_\tau$ with $\tau=2$ when $|\mathcal{T}|=10$. The blue and orange lines represent the fitting results of Gaussian-assumption-based model and GMM, respectively. Obviously, the uncertainties do not follow Gaussian distribution but can be well fitted by GMM.}
	\label{fig_sample}
	\vspace{-4mm}
\end{figure}

Based the whole-year regD signals from PJM in 2020 \cite{dispatch2019pjm}, the uncertainties $\overline u_\tau$ and $\underline u_\tau$ do not follow Gaussian distribution, as shown in Fig. \ref{fig_sample}. Moreover, these uncertainties are on the LHS in  (\ref{eqn_thermal_r}) because they are multiplied with a decision variable $R^\text{ha}$. Thus, Eq. (\ref{eqn_thermal_r}) is still hard to deal with. 




\subsection{Mixture-model-based convexification approach}
To address the intractability caused by the non-Gaussian LHS uncertainties in  (\ref{eqn_thermal_r}), a mixture-model-based convexification approach is proposed.
We first introduce GMM to fit the original non-Gaussian uncertainties with multiple Gaussian distributions. Then, we reformulate  (\ref{eqn_thermal_r}) into deterministic non-convex forms. {Finally, we propose tractable SOCP approximations for these non-convex constraints.}

\subsubsection{Introduction of GMM}
GMM can approximate the distribution of non-Gaussian variable $\bm \omega$ with multiple Gaussian distributions \cite{goodfellow2016deep}:
\begin{align}
p^\text{NG}(\bm \omega)=\sum_{j \in \mathcal{J}} \pi_{j} p\left(\bm \xi_{j}|\bm \mu_{j}, \bm \Sigma_{j}\right), \label{eqn_GMM}
\end{align}
where $p^\text{NG}(\bm \omega)$ is the probability density function (PDF) of $\bm \omega$; $j \in \mathcal{J}$ is the index of the Gaussian component; $\pi_{j}$ is the weight of component $j$, and $\sum_{j \in \mathcal{J}} \pi_{j}=1$; $p\left(\bm \xi_{j}|\bm \mu_{j}, \bm \Sigma_{j}\right)$ represents the PDF of a Gaussian uncertainty $\bm \xi_{j}$; $\bm \mu_{j}$ and $\bm \Sigma_{j}$ are the expectation and covariance of $\bm \xi_{j}$, respectively. 
{Based on the historical samples of $\bm \omega$, the three parameters $\pi_{j}$, $\bm \mu_{j}$, and $\Sigma_{j}$ can be estimated based on the Expectation Maximization algorithm \cite{8017474,9122389,9535415}.} Fig. \ref{fig_sample} provides an example to demonstrate the excellent fitting power of GMM. 
  
\subsubsection{Deterministic reformulations of chance constraints}
The generic form of the chance constraints in  (\ref{eqn_thermal_r}) can be expressed as follows:
\begin{align}
\mathbb{P}\left(\bm \alpha(\bm x)^\intercal \bm \omega \leq \beta(\bm x)\right) \geq 1-\epsilon. \label{eqn_cc_generic}
\end{align}
The detail expressions of $\bm \alpha(\bm x)$ and $\beta(\bm x)$ for each constraint is given in Appendix \ref{app_1}.
To reformulate  (\ref{eqn_cc_generic}), we introduce the following \textbf{Lemma}.
\begin{lemma} \label{lemma_1}
If the PDF of the uncertainty $\bm \omega$ is approximated by GMM, i.e., Eq. (\ref{eqn_GMM}), then we have \cite{hu2018chance}
\begin{align}
\mathbb{P} \left( \bm \alpha(\bm x)^\intercal \bm \omega \leq \beta(\bm x)\right)=\sum_{j \in \mathcal{J}}\pi_j\mathbb{P} \left( \bm \alpha(\bm x)^\intercal \bm \xi_{j} \leq \beta(\bm x)\right).
\end{align}
\end{lemma}
By introducing an auxiliary variable $y_{j}$ for each Gaussian component, Eq. (\ref{eqn_cc_generic}) can be converted into:
\begin{align}
&\mathbb{P} \left( \bm \alpha(\bm x)^\intercal \bm \xi_{j} \leq \beta(\bm x)\right) \geq y_j, \quad \forall j \in \mathcal{J}, \label{eqn_cc_r1}\\
&\sum_{j \in \mathcal{J}}\pi_{j}y_{j} \geq 1-\epsilon. \label{eqn_probability}
\end{align}
Note $y_{j}$ should be no less than 0.5 due to the definition of the CCP. 
Since the uncertainty $\bm \xi_{j}$ follows Gaussian distribution, Eq. (\ref{eqn_cc_r1}) can be further reformulated into the following deterministic non-convex form:
\begin{align}
\Phi^{-1}(y_j)\sqrt{\bm \alpha(\bm x)^\intercal \bm \Sigma_j \bm \alpha(\bm x)} + \alpha(\bm x)^\intercal \bm \mu_j \leq \beta(\bm x), \forall j \in \mathcal{J}, \label{eqn_cc_r2}
\end{align}
where $\Phi^{-1}(\cdot)$ is the inverse of the cumulative distribution function of the stand normal distribution.

\begin{remark}
Eq. (\ref{eqn_cc_r2}) is still intractable. On the one hand, the term $\Phi^{-1}(y_j)\sqrt{\bm \alpha(\bm x)^\intercal \bm \Sigma_j \bm \alpha(\bm x)}$ is non-convex because both $y_j$ and $\bm \alpha(\bm x)$ are variables. On the other hand, function $\Phi^{-1}(\cdot)$ has no analytical formula.
\end{remark}

\subsubsection{Convexification for deterministic reformulations}
\begin{figure}
	\centering
	\vspace{-4mm}
	\includegraphics[width=0.9\columnwidth]{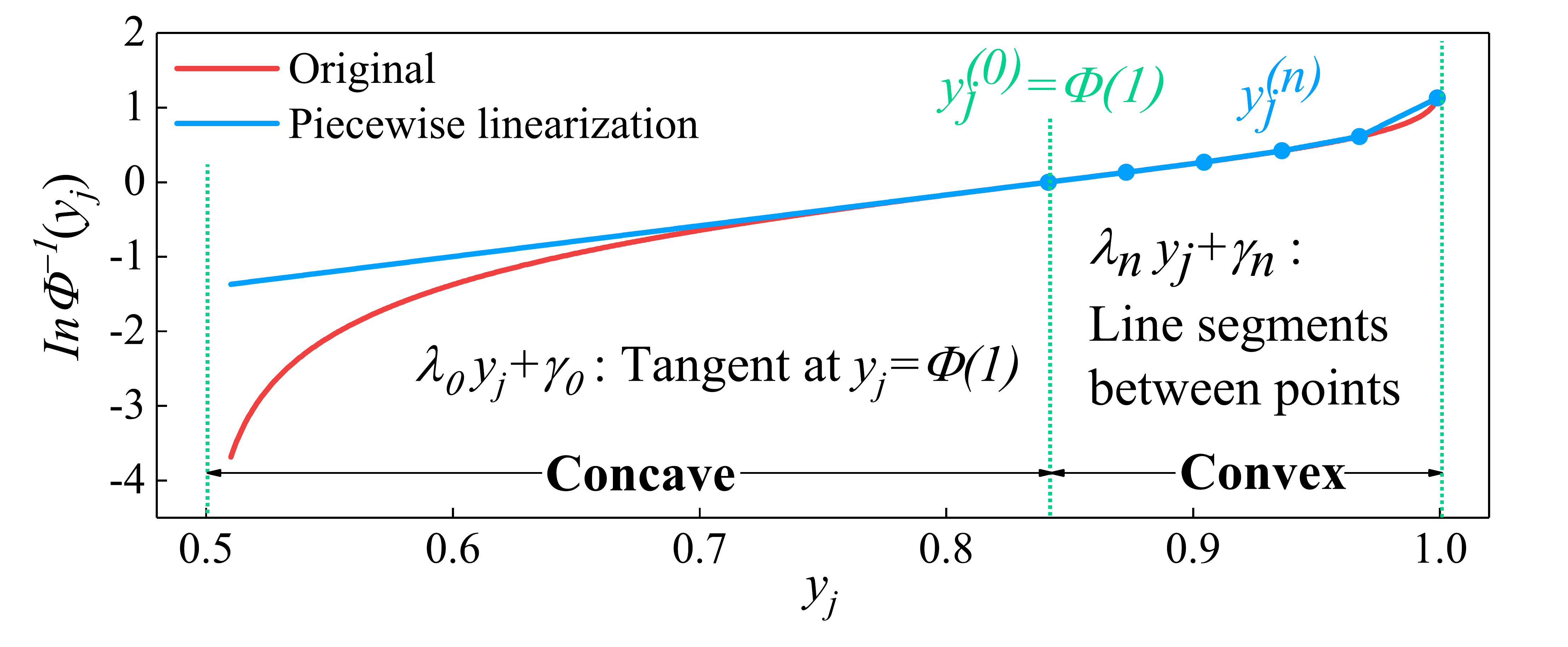}\vspace{-4mm}
	\caption{Graphs of function $\ln{\Phi^{-1}(y_j)}$ (red line) and the proposed inner approximation based on piecewise linear function $\max_{n \in \mathcal{N}}{\{\lambda_n y_j + \gamma_n}\}$ (blue line). The blue dots, i.e., $(y_j^{(n)}, \ln{\Phi^{-1}(y_j^{(n)})})$, are points on function $\ln{\Phi^{-1}(y_j)}$. The first line, i.e., $\lambda_0 y_j + \gamma_0$ is the tangent of $\ln{\Phi^{-1}(y_j)}$ at $y_j = \Phi(1)$. The rest lines are segments between two adjacent blue dots. Note the piecewise linear function $\max_{n \in \mathcal{N}}{\{\lambda_n y_j + \gamma_n}\}$ is convex.
	}
	\label{fig_PWL}
	 		\vspace{-4mm}
\end{figure}

We propose a convexification method to approximately reformulate the deterministic reformulation (\ref{eqn_cc_r2}) into a tractable SOCP form. 
Firstly, observing that all elements of $\bm \alpha(\bm x)$ in our problem are always nonegative (See Appendix \ref{app_1}), we can re-express each element of $\bm \alpha(\bm x)$ in an exponential manner:
\begin{align}
\alpha_k = e^{\rho_k}, \quad \forall k \in \mathcal{K}, \label{eqn_exp_replace}
\end{align}
where $\alpha_k$ is the $k$-th element of $\bm \alpha(\bm x)$; $\rho_k$ is an auxiliary variable; $\mathcal{K}$ is the corresponding index set. Since $y_j \geq 0.5$, function $\Phi^{-1}(y_j)$ is nonegative, so it can be also expressed as an exponential form, i.e., $\Phi^{-1}(y_j) = e^{\ln{\Phi^{-1}(y_j)}}$. According to the statistics for the whole-year regulation signals in 2020 \cite{dispatch2019pjm}, the off-diagonal elements in the covariance matrix $\bm \Sigma_j$ is close to zero and much smaller than the diagonal ones.
Thus, we can use zero to replace the off-diagonal elements of $\bm \Sigma_j$. Then, by substituting the above exponential expressions, the first term in (\ref{eqn_cc_r2}) can be converted into a L2-norm form:
\begin{align}
\Phi^{-1}(y_j)&\sqrt{\bm \alpha(\bm x)^\intercal \bm \Sigma_j \bm \alpha(\bm x)} = \sqrt{\sum_{k}(\sigma_{k}e^{\rho_k + \ln{\Phi^{-1}(y_j)}}})^2 \notag\\
&= \Vert \sigma_{k}e^{\rho_k + \ln{\Phi^{-1}(y_j)}}, \quad \forall k \in \mathcal{K} \Vert_2, \quad \forall j \in \mathcal{J}, \label{eqn_norm_1}
\end{align}
where $\sigma_{k}$ is the standard deviation of the $k$-th element of $\bm \omega$, and its value is always non-negative. 
Since both the exponential and L2-norm functions are convex and element-wise monotonically increasing, the L2-norm term in (\ref{eqn_norm_1}) is convex if its power, $\rho_k + \ln{\Phi^{-1}(y_j)}$, is convex according to convex condition for composite functions \cite{boyd2004convex}. Unluckily, this power is non-convex in its domain. Nevertheless, we can find a piecewise-linearization-based safe approximation to convexify this non-convex power based on the following \textbf{Proposition}.
\begin{proposition} \label{prop_1}
In the domain of function $\ln{\Phi^{-1}(y_j)}$, i.e., $y_j>0.5$, we have
\begin{align}
\ln{\Phi^{-1}(y_j)} \leq \max_{n \in \mathcal{N}}{\{\lambda_n y_j + \gamma_n}\}, \quad \forall j \in \mathcal{J}, \label{eqn_prop_1}
\end{align}
where $\mathcal{N}=\{0, 1, \cdots, N\}$ is the index set of lines. Symbol $\lambda_n y_j + \gamma_n, \forall n \in \mathcal{N}$ represents different lines constructed by piecewise linearization (shown in Fig. \ref{fig_PWL}), i.e., $\lambda_0 y_j + \gamma_0$ is the tangent at $y_j = \Phi(1)$, while the rest lines, i.e., $\lambda_n y_j + \gamma_n, \forall n \in \mathcal{N}/\{0\}$, are line segments by connecting two points on function $\ln{\Phi^{-1}(y_j)}$ in sequence.  Note the piecewise linear function $\max_{n \in \mathcal{N}}{\{\lambda_n y_j + \gamma_n}\}$ is convex because its epigraph is a convex polyhedron.
\end{proposition}
\emph{Proof:} See Appendix \ref{app_2}. 

Observing that both the exponential and L2-norm functions are monotonically increasing, the following inequality holds according to \textbf{Proposition} \ref{prop_1}: 
\begin{align}
\Vert \sigma_{k}&e^{\rho_k + \ln{\Phi^{-1}(y_j)}}, \forall k \in \mathcal{K} \Vert_2 \notag\\
&\leq \Vert \sigma_{k}e^{\rho_k + \max_{n \in \mathcal{N}}{\{\lambda_n y_j + \gamma_n}\}},\forall k \in \mathcal{K} \Vert_2, \forall j \in \mathcal{J}. \label{eqn_approximation_UB}
\end{align}
Based on (\ref{eqn_norm_1}) and (\ref{eqn_approximation_UB}), the deterministic non-convex reformulation (\ref{eqn_cc_r2}) can be convexified as the following SOCP form:
\begin{align}
\Vert \sigma_{k}e^{\rho_k + \max_{n \in \mathcal{N}}{\{\lambda_n y_j + \gamma_n}\}},\forall k \in \mathcal{K} \Vert_2 + &\bm \alpha(\bm x)^\intercal \bm \mu_{j} \leq \beta(\bm x), \notag \\
& \quad \quad \forall j \in \mathcal{J}, \label{eqn_r_socp_UB_0}
\end{align}
which is also equivalent to
\begin{align}
\Vert \sigma_{k}e^{\rho_k + \lambda_n y_j + \gamma_n},\forall k \in \mathcal{K} \Vert_2 + &\bm \alpha(\bm x)^\intercal \bm \mu_{j} \leq \beta(\bm x), \notag \\
&\forall n \in \mathcal{N}, \quad \forall j \in \mathcal{J}. \label{eqn_r_socp_UB}
\end{align}

\begin{remark}
The inner approximation used in (\ref{eqn_approximation_UB})-(\ref{eqn_r_socp_UB}) introduces additional conservativeness due to the approximation error. According to Fig. \ref{fig_PWL} and \textbf{Proposition} \ref{prop_1}, this error is small {when $y_j>\Phi(1)$} because we introduce multiple line segments for approximation. Only when $y_j$ is close to 0.5, this error becomes significant. Nevertheless, Eq. (\ref{eqn_probability}) requires that the weighted average of $y_j$ should be no smaller than $1-\epsilon$. The risk parameter $\epsilon$ is usually small, so $y_j$ will be large. Therefore, the additional conservativeness introduced by the approximation error is commonly insignificant. This ensures the optimality performance of the proposed method.
\end{remark}

In the previous convexification, we introduce another non-convex constraint, i.e., Eq. (\ref{eqn_exp_replace}). Nevertheless, the vector $\bm \alpha(\bm x)$ only contains one single variable $R_t^\text{ha}$ (See Appendix \ref{app_1}), leading to only one non-convex constraint, i.e., $R_t^\text{ha} = e^{\rho_2}$. We can also use piecewise linearization to reformulate this non-convex constraint into:
\begin{align}
R_t^\text{ha}=\max_{n \in \mathcal{N}^{R}}{\lambda_n^{R} \rho_2 + \gamma_n^{R}}, \label{eqn_PWL_R}
\end{align}
where $\mathcal{N}^{R}$ is the index set of lines for piecewise linearization. The $n$-th
line, i.e., $\lambda_n^{R} \rho_2 + \gamma_n^{R}$, is constructed by connecting the $n$-th and $(n+1)$-th points on function $e^{\rho_2}$. Unlike (\ref{eqn_r_socp_UB_0}), the maximum operator appears on the RHS of ``$=$" in (\ref{eqn_PWL_R}). Thus, we need to employ the Big-M method with auxiliary binary variable $z_n$ to reformulate (\ref{eqn_PWL_R}) into a solvable form:
\begin{align}
\begin{cases}
&R_t^\text{ha} \geq \lambda_n^{R} \rho_2 + \gamma_n^{R}, \quad \forall n \in \mathcal{N}^{R}, \\
&R_t^\text{ha} \leq \lambda_n^{R} \rho_2 + \gamma_n^{R} + M(1-z_n), \quad \forall n \in \mathcal{N}^{R}, \\
& \bm {1}^{\intercal} \bm z = 1, \quad \bm z \in \{0,1\}^{|\mathcal{N}^{R}|},
\end{cases} \label{eqn_PWL_R2}
\end{align}
where $M$ is a big number and $|\mathcal{N}^{R}|$ is the length of set $\mathcal{N}^{R}$. 

\begin{remark}
The auxiliary binary variable number, i.e., $|\mathcal{N}^{R}|$, can be very small because only one single variable $R_t^\text{ha}$ needs to be linearized based on  (\ref{eqn_PWL_R2}), which guarantees computational tractability.
\end{remark}

Finally, since  (\ref{eqn_obj}) is linear, the expectation of the total cost, $\mathbb E (EC_t)$, can be calculated by: 
\begin{align}
    \mathbb E (EC_t) = &\sum_{i \in \mathcal{I}} \eta_t(p_{t,i}^\text{ha}- R_{t,i}^\text{ha} s_t^\text{avg})\Delta t \notag \\
    &\quad \quad - \sum_{i \in \mathcal{I}} (r_t^\text{rc} + r_t^{m}m_t^\text{avg})R_{t,i}^\text{ha}, \label{eqn_obj_r}
\end{align}
where $s_t^\text{avg}=\mathbb{E}(s_{t,l})=\frac{\sum_{l \in \mathcal{L}}s_{t,l}}{|\mathcal{L}|}$ and $m_t^\text{avg}=\mathbb{E}(m_t)$. Then, \textbf{P1} can be reformulated into a mixed-integer SOCP problem:
\begin{align} 
&\min \quad \text{Eq. (\ref{eqn_obj_r})} \tag{$\textbf{P2}$},\\
&\begin{array}{r@{\ }r@{}l@{\ }l}
\text{s.t.:} &&\text{Eqs.} &\text{ (\ref{eqn_r_power})-(\ref{eqn_R_bound}), \{(\ref{eqn_probability}), (\ref{eqn_r_socp_UB})\}$_{(\ref{eqn_thermal_r})}$, and (\ref{eqn_PWL_R2})},
\end{array} \notag
\end{align}
where \{(\ref{eqn_probability}), (\ref{eqn_r_socp_UB})\}$_{(\ref{eqn_thermal_r})}$ represents that each chance constraint in  (\ref{eqn_thermal_r}) is reformulated into  (\ref{eqn_probability}) and (\ref{eqn_r_socp_UB}) in \textbf{P2}.

\section{Case study} \label{sec_case}
\subsection{System Configuration}
\begin{figure}
	\subfigbottomskip=-6pt
	\subfigcapskip=-4pt
	\centering
	\subfigure[]{\includegraphics[width=0.9\columnwidth]{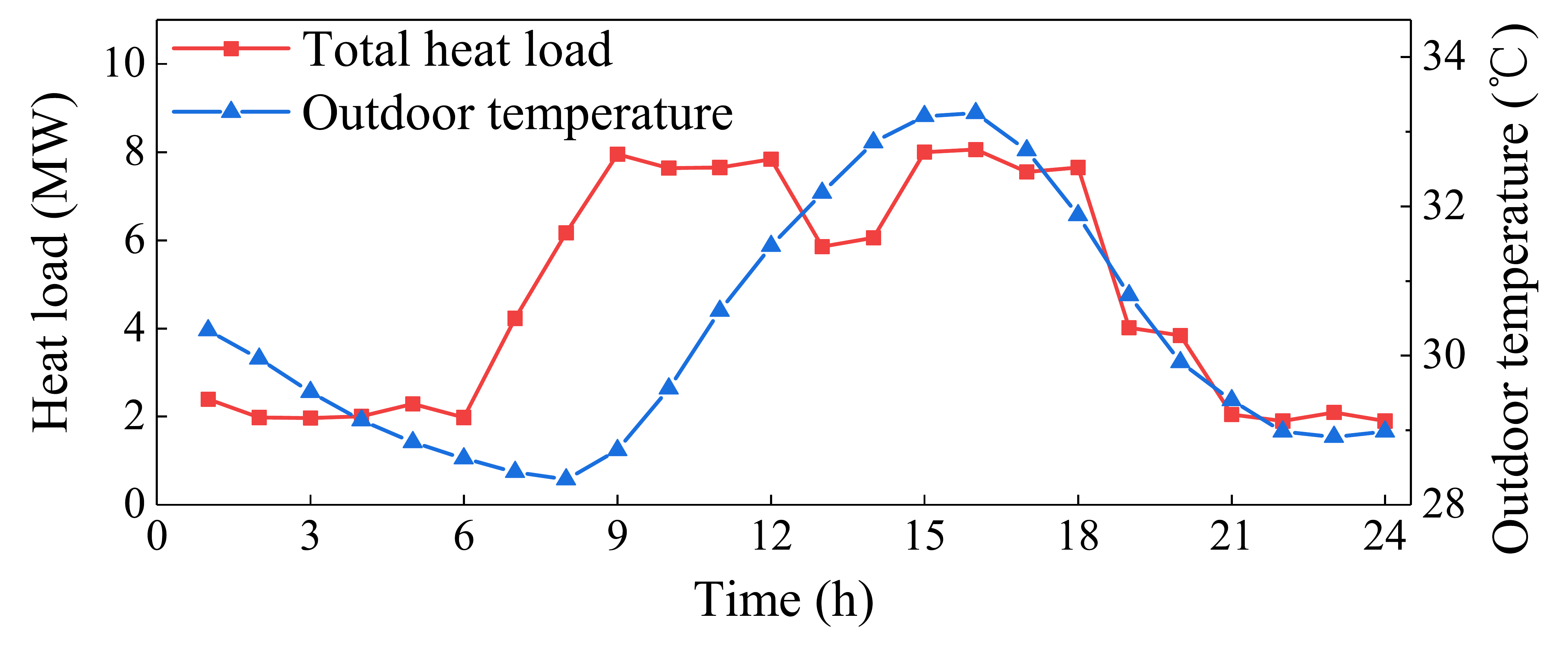}}
	\subfigure[]{\includegraphics[width=0.9\columnwidth]{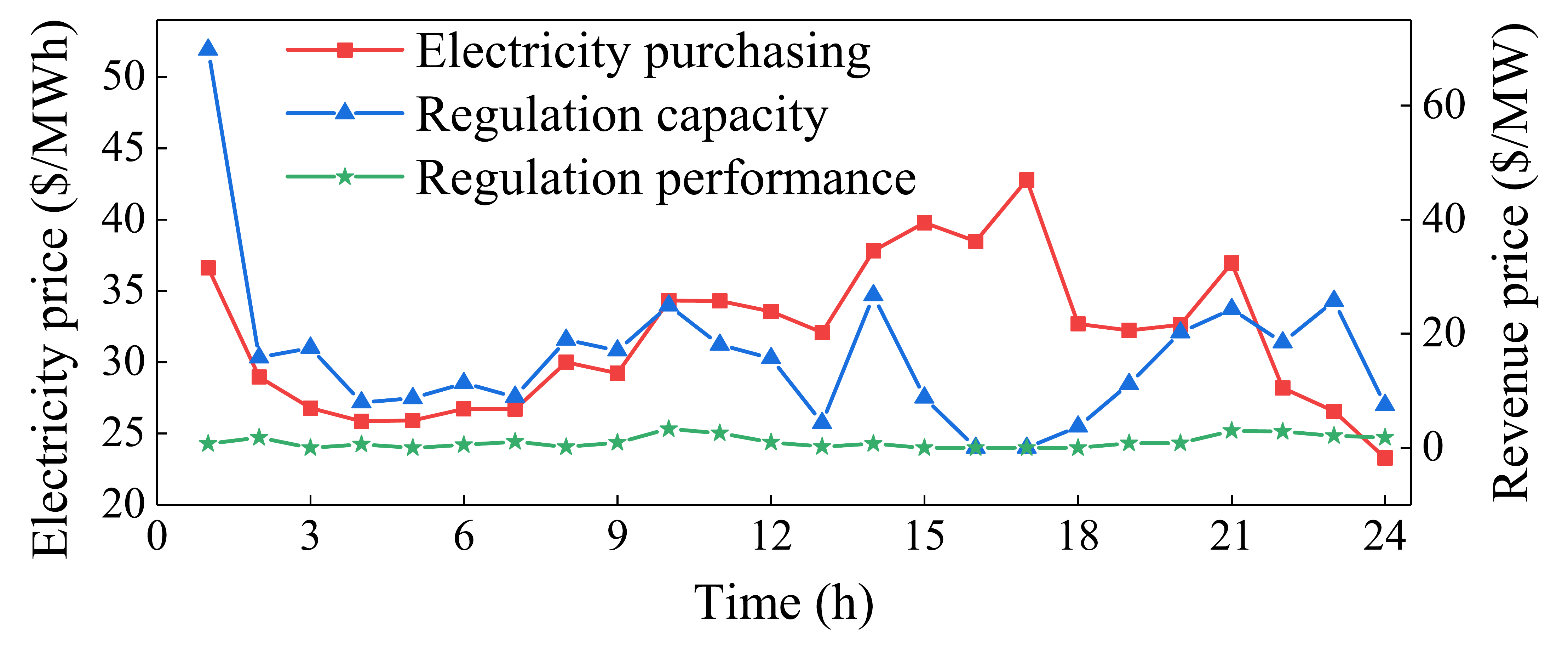}}
	\vspace{1mm}
	\caption{(a) Heat load $h_t$ and outdoor temperature $\theta^\text{out}$ and (b) electricity purchase price $e_t$, regulation capacity revenue price $r_t^\text{rc}$, and regulation performance revenue price $r_t^\text{m}$.}
	\label{fig_load}
\end{figure}

We validate the proposed method based on one large-capacity HVAC system.
The daily heat loads, outdoor temperature, and unit prices for electricity purchasing and regulation revenue are demonstrated in Fig. \ref{fig_load}.
Other parameters are listed in Table \ref{tab_parameter}. 
We collect the whole-year regD signals from PJM in 2020 \cite{dispatch2019pjm} as historical data. Based on these signal data, we construct the samples of the uncertain parameters used in \textbf{P2}, including $m_t$ in (\ref{eqn_millage}), $\overline u_\tau$ and $\underline u_\tau$ in (\ref{eqn_u_extreme}). The whole dataset has been uploaded in \cite{samples2021}. 

All simulations are implemented based on an Intel(R) Core(TM) 8700 3.20GHz CPU with 16GB memory. The corresponding optimization problem is built by CVXPY and solved by MOSEK. 

\begin{table}
	\small
	\centering
	\vspace{-4mm}
	\caption{Parameters in case study}
	\vspace{-2mm}
	\begin{threeparttable} 
	\begin{tabular}{cccc}
		\hline
		\rule{0pt}{11pt}		
		Parameters & Value &Parameters & Value\\
		\hline
		\rule{0pt}{10pt}
		$C_i$ &  1.75 MWh/{\textcelsius} & $\theta^\text{min}$& 22{\textcelsius}\\
		$g_i$ &  0.2MW/{\textcelsius} & $\theta^\text{max}$& 28{\textcelsius} \\
		$COP_i$ &  5 &  $P_i^\text{max}$ & 2MW \\
		$|\mathcal{N}|^*$ & 10 & $P_i^\text{min}$ & 0 \\
		${|\mathcal{N}^{R}|}^*$ & 50 & $|\mathcal{T}|^*$ & 10\\
		$|\mathcal{J}|^*$ & 3 & & \\
		\hline
	\end{tabular}\label{tab_parameter}
	\begin{tablenotes}
	\footnotesize
	\item[*] Symbol $|\cdot|$ denotes the length of the set.
	\end{tablenotes}
	\end{threeparttable}
	\vspace{-2mm}
\end{table}

\subsection{Benchmarks}
To demonstrate the superiority of the proposed approach, three benchmarks are introduced: 
\begin{enumerate}
	\item \textbf{B1}: The Gaussian-assumption-based CCP used in \cite{8017474,9122389,9535415};
	\item \textbf{B2}: The moment-based DRCC method used in \cite{8301597};
	\item \textbf{B3}: The Wasserstein-distance-based DRCC method used in \cite{9347818}. In this approach, the proposed temporal compression method is also applied; otherwise the computational burden will be too heavy and out of memory issue will occur. A total of 1000 samples are used for constructing the ambiguity set\footnote{The optimality of \textbf{B3} can be improved by increasing the sample number for constructing the ambiguity set \cite{esfahani2018data}. However, the computational burden is also proportional to this sample number. According to our test, if this sample number is larger than 1000 (e.g. 1500), out of memory issue occurs.}.
\end{enumerate}

\subsection{Model comparison}
\subsubsection{Optimality, feasibility, and computational efficiency}
\begin{figure}
	\subfigbottomskip=-6pt
	\subfigcapskip=-4pt
	\centering
	\subfigure[]{\includegraphics[width=0.9\columnwidth]{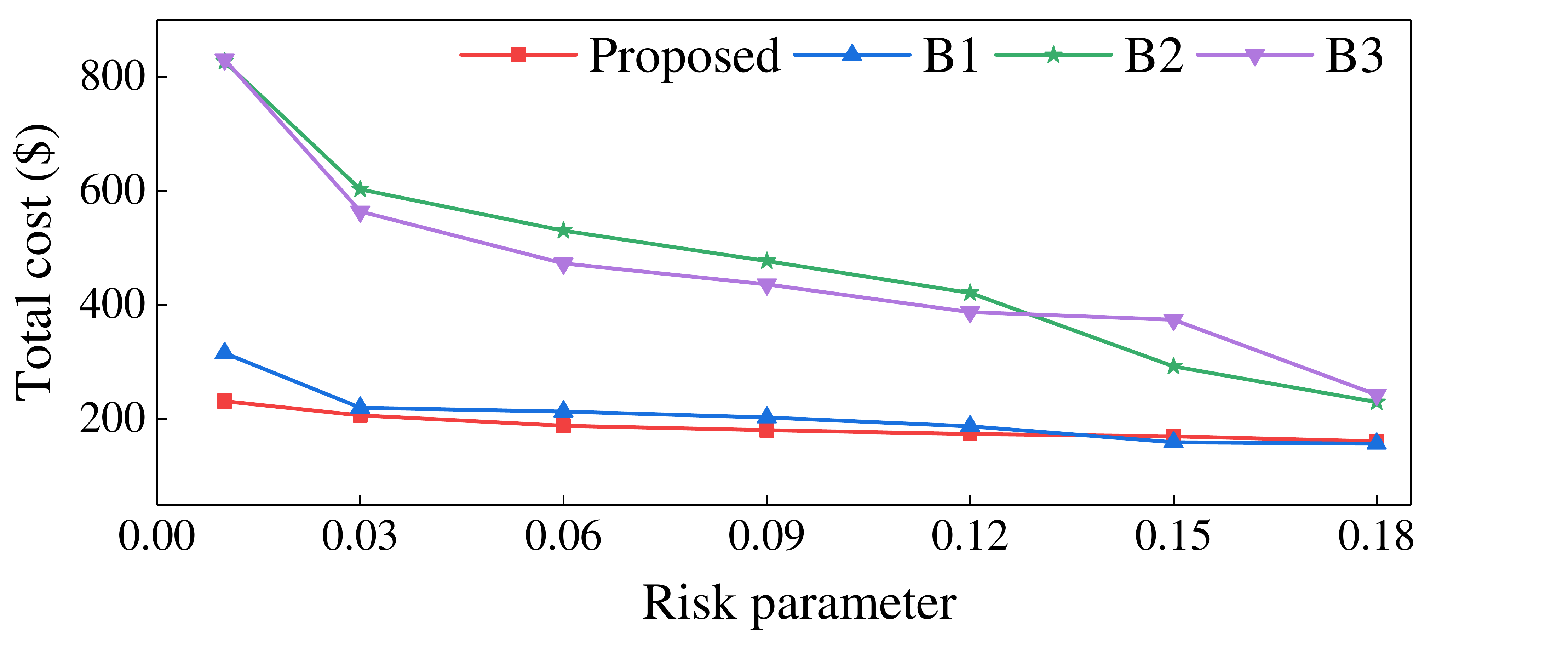}}
	\subfigure[]{\includegraphics[width=0.9\columnwidth]{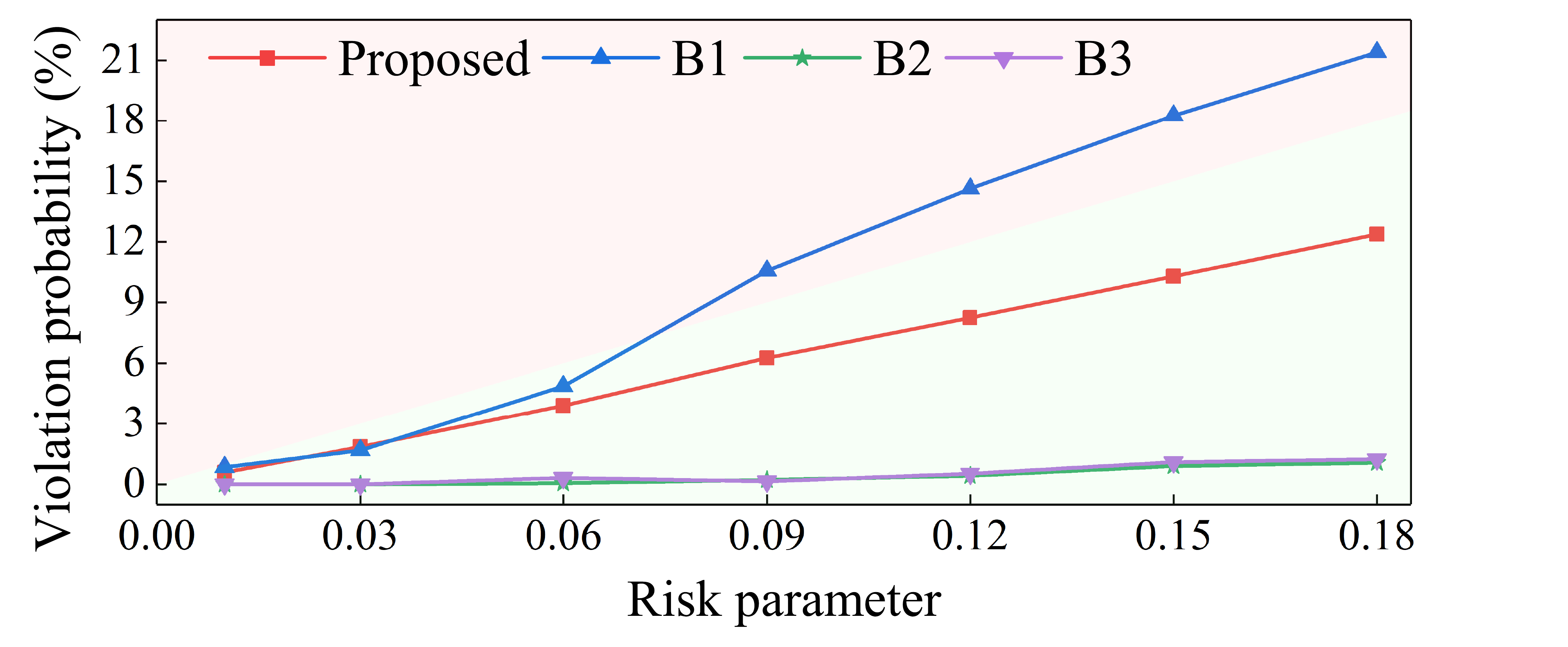}}
	\subfigure[]{\includegraphics[width=0.9\columnwidth]{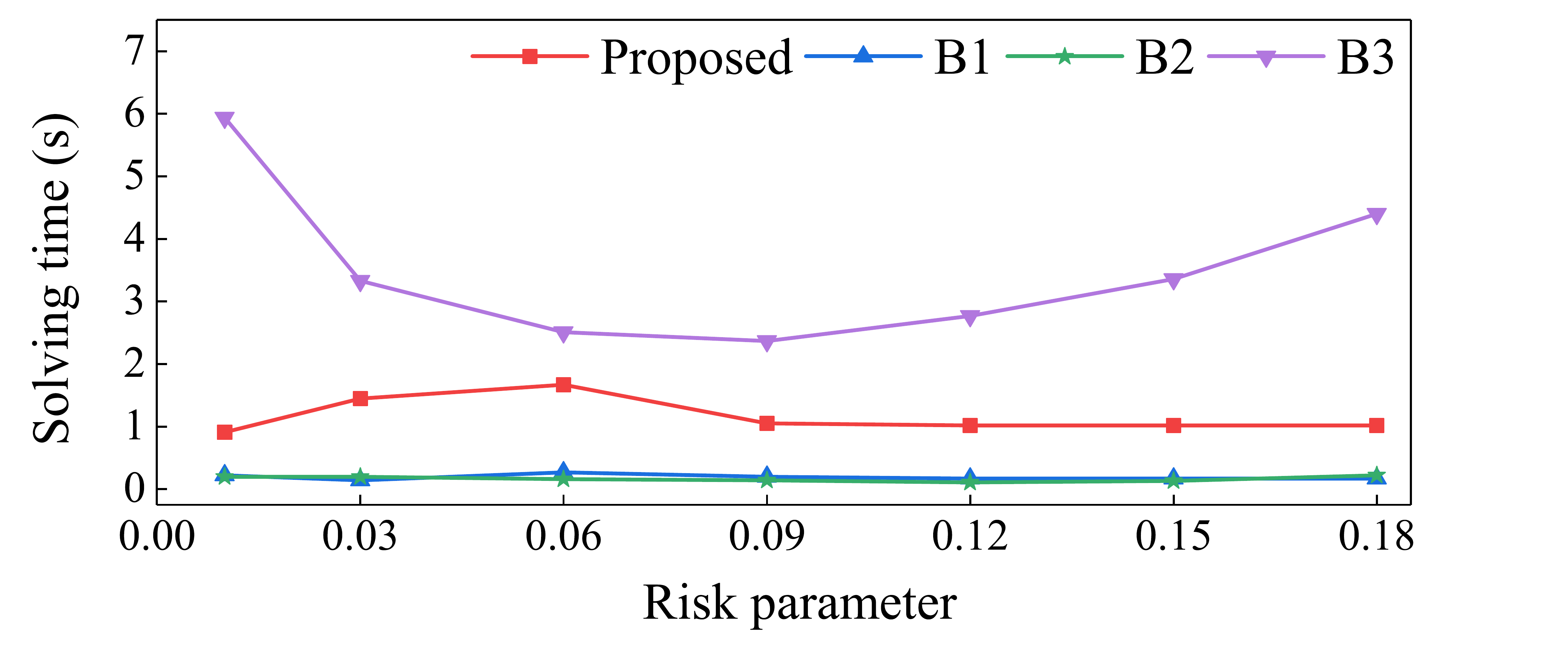}}
	\vspace{1mm}
	\caption{Results of (a) whole-day total cost, i.e., $\sum_{t=1}^{24}EC_t$, (b) maximum violation probability, and (c) solving time obtained by different models. In (b), the green and red regions represent the safe (i.e. the maximum violation probability is smaller than the risk parameter) and unsafe regions (i.e. the maximum violation probability is larger than the risk parameter).  }
	\label{fig_comparison}
	\vspace{-4mm}
\end{figure}

Fig. \ref{fig_comparison} shows the results of the whole-day total costs (i.e. $\sum_{t=1}^{24}EC_t$), solving times, and maximum probability violation under different risk parameters. In all cases, the total cost of the proposed method is almost the lowest, while its maximum violation probability always stays in the allowable range. Although the Gaussian-assumption-based method \textbf{B1} can achieve comparable optimality to the proposed one, it may derive infeasible solutions, i.e., the maximum violation probability is larger than the given risk parameter. This is because the original uncertainties are non-Gaussian distributed.
The rest two DRCC methods, i.e. \textbf{B2} and \textbf{B3}, can always ensure the feasibility of solutions. However, their total costs are much higher than that of the proposed method because they need to ensure the feasibility for all the possible distributions in their ambiguity sets, including some distributions that are quite different from the actual one. 

In \textbf{B1} and \textbf{B2}, each chance constraint is reformulated into only one single SOCP constraint but introduces no additional constraint, while some additional constraints, such as  (\ref{eqn_PWL_R2}), are necessary for the proposed method. 
Thus, the solving times of \textbf{B1} and \textbf{B2} are lower than that of the proposed model. Nevertheless, the solving time of the proposed method is only around 1s, which is also acceptable in practice. In \textbf{B3}, numerous additional constraints have to be introduced \cite{esfahani2018data}, so its computational performance is the worst among all methods. These results confirm the great optimality and feasibility performance of the proposed method.

\subsubsection{Hour-ahead regulation capacity offers}
\begin{figure}
	\centering
	\includegraphics[width=0.9\columnwidth]{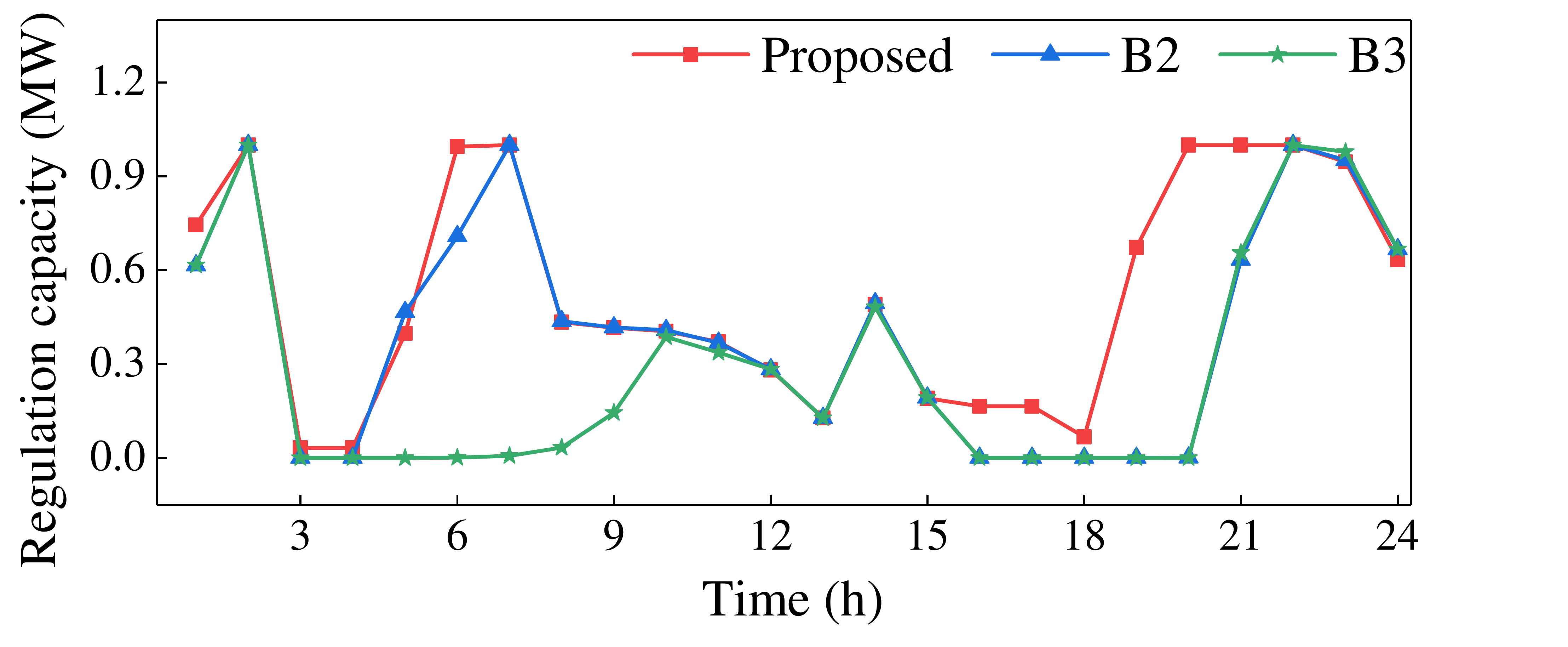}\vspace{-4mm}
	\caption{Results of hour-ahead regulation capacity offers under $\epsilon=0.15$.  
	}
	\label{fig_R_ha}
	 		\vspace{-2mm}
\end{figure}

Fig. \ref{fig_R_ha} summarizes the hour-ahead regulation capacity offers obtained by different methods under $\epsilon=0.15$. Note the results of \textbf{B1} are not listed here because \textbf{B1} can not ensure the feasibility of solutions. In all time, the hour-ahead regulation capacity offer of the proposed method is much larger compared to the rest models. As aforementioned, \textbf{B2} and \textbf{B3} have to satisfy constraints for all distributions in their ambiguity sets, so they are more conservative. As a result, large margins need to be reserved for the uncertainties in both \textbf{B2} and \textbf{B3}, which shrinks the potential regulation capacity. This result validates the 
better optimality of the proposed method.

\subsection{Sensitivity analysis} \label{sec_sensitivity}
\subsubsection{Time duration number $|\mathcal{T}|$ for splitting the operating hour}
\begin{figure}
	\subfigbottomskip=-6pt
	\subfigcapskip=-4pt
		\vspace{-4mm}
	\centering
	\subfigure[]{\includegraphics[width=0.9\columnwidth]{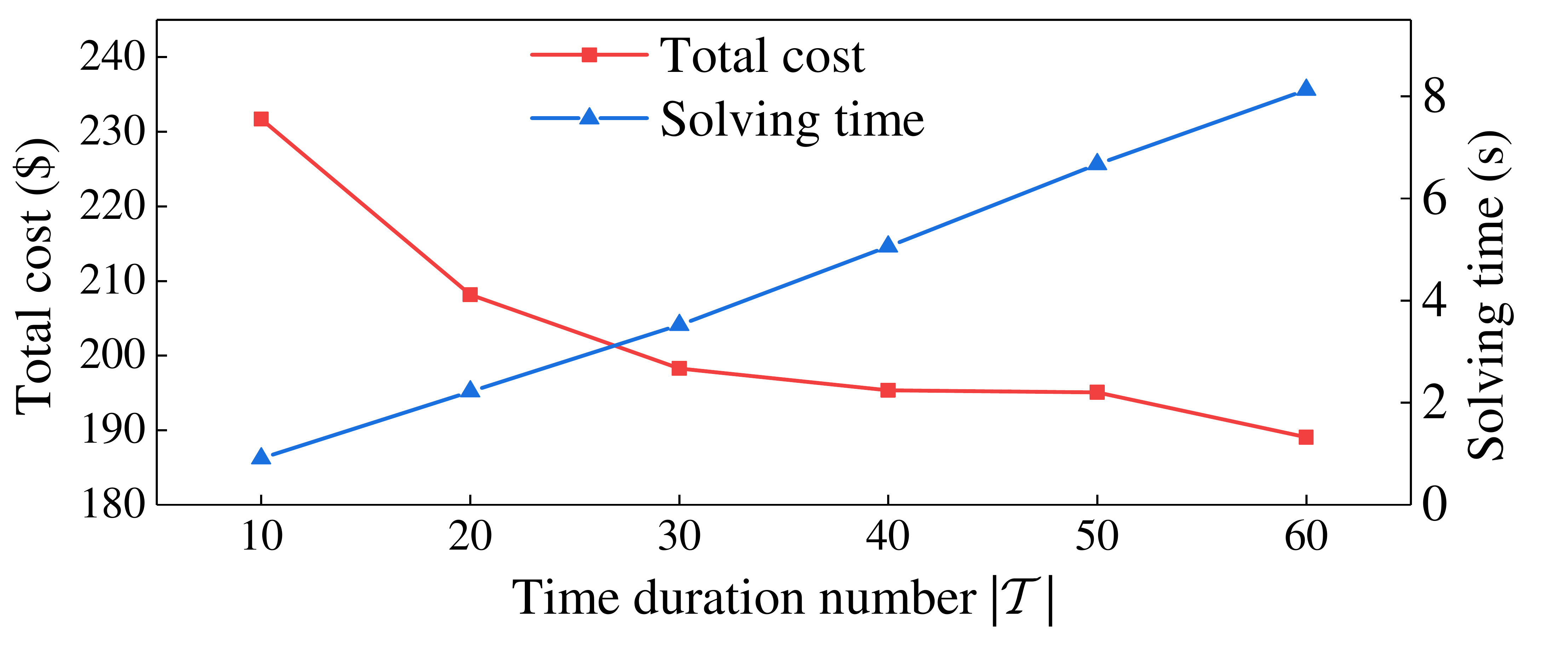}}
	\subfigure[]{\includegraphics[width=0.9\columnwidth]{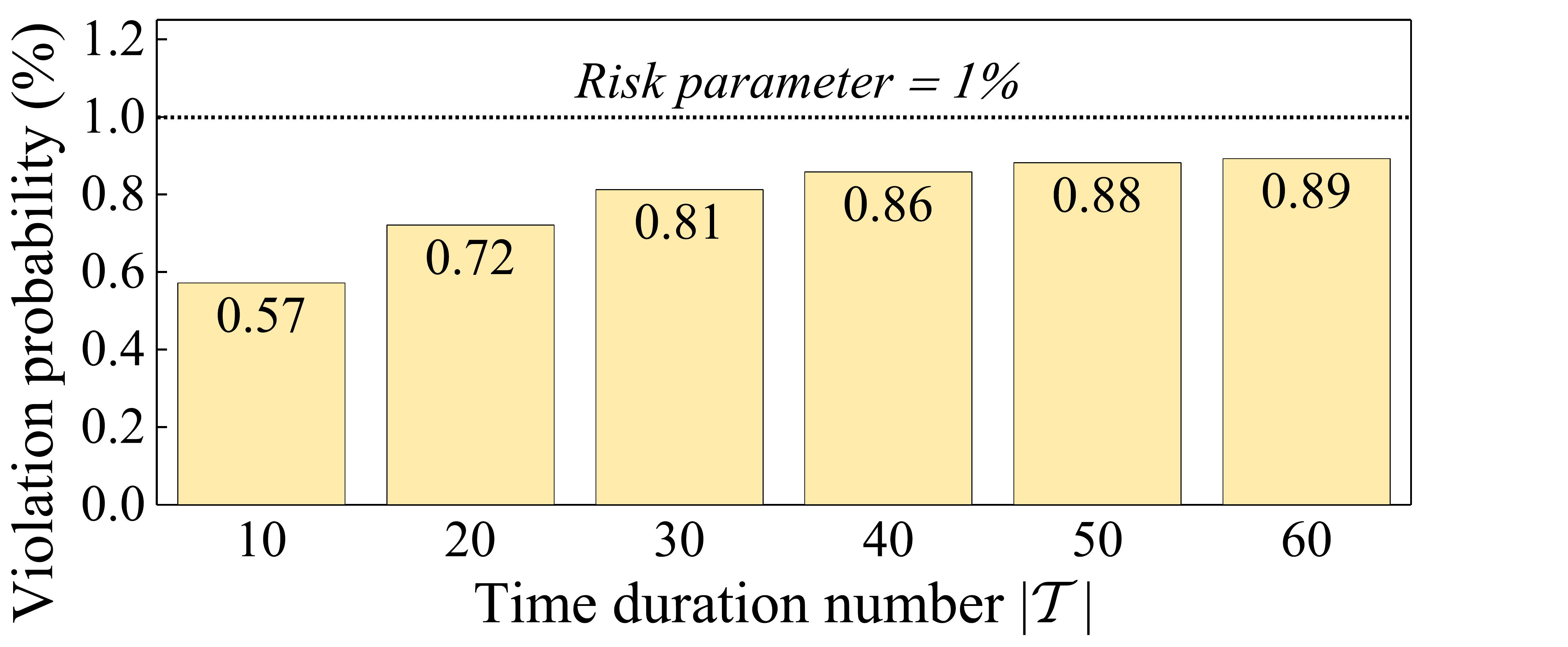}}
	\vspace{1mm}
	\caption{Results of (a) daily total cost and solving time, and (b) maximum violation probability under different time duration numbers, i.e., $|\mathcal{T}|$, for splitting the operating hour with $\epsilon=0.01$. }
	\label{fig_duration}
	\vspace{-4mm}
\end{figure}

We implement a case study with different $|\mathcal{T}|$ to investigate its effects, and the corresponding results are shown in Fig. \ref{fig_duration}. The risk parameter is 0.01, while the line number for the piecewise linearization $|\mathcal{N}|^R$ is set as 50. With the growth of $|\mathcal{T}|$, the total cost decreases. According to \textbf{Proposition} \ref{prop_0}, if one time duration $\Delta \tau$ is split into smaller ones, the additional conservativeness introduced by the inner approximation used in (\ref{eqn_approximation}) can be further reduced. That is to say, increasing $|\mathcal{T}|$ can improve the optimality of the proposed method. Conversely, the solving time becomes larger because more constraints are introduced according to (\ref{eqn_r_socp_UB}). Because increasing $|\mathcal{T}|$ decreases the conservativeness, the maximum violation probability also grows with the increase of $|\mathcal{T}|$. Nevetheless, its value always maintains in the safe region, which demonstrates the great feasibility of the proposed method.

\subsubsection{Line number $|\mathcal{N}|^R$ in (\ref{eqn_PWL_R})} 
\begin{figure}
	\subfigbottomskip=-6pt
	\subfigcapskip=-4pt
	\centering
	\subfigure[]{\includegraphics[width=0.9\columnwidth]{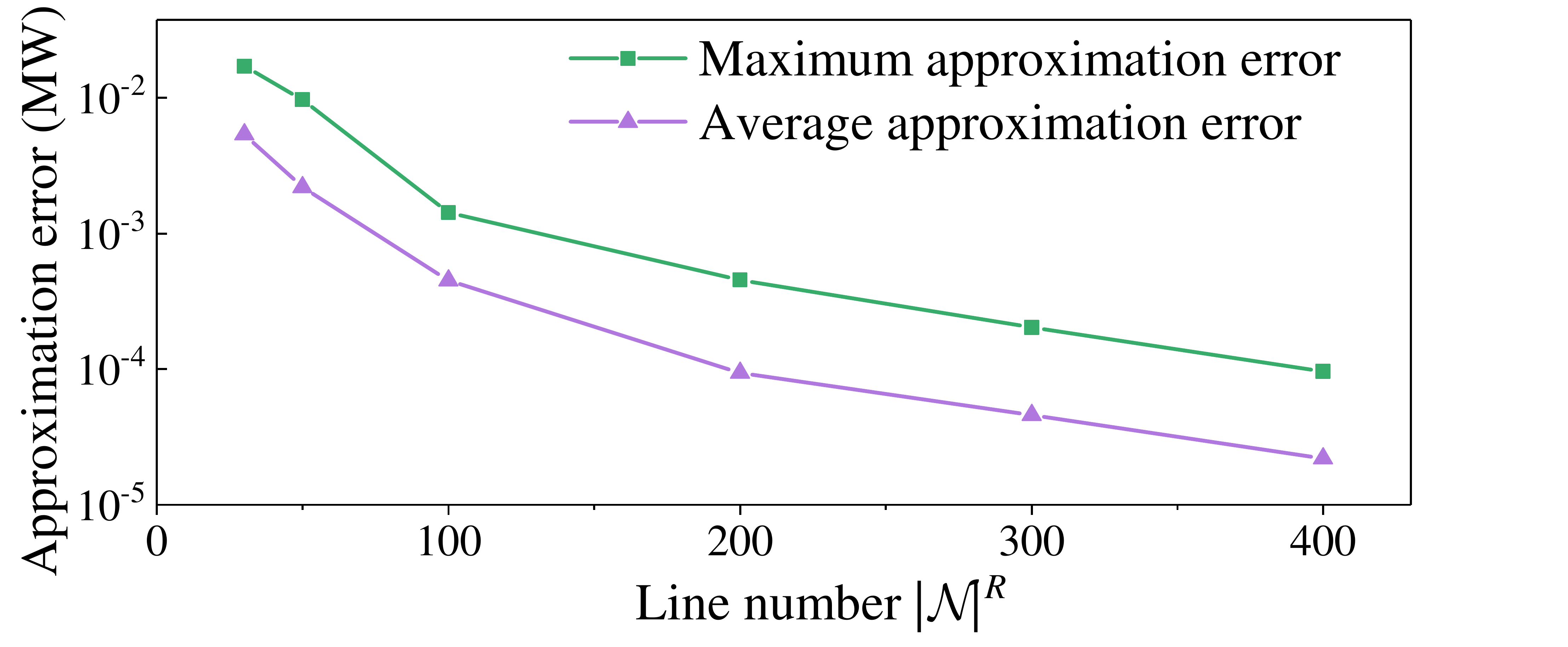}}
	\subfigure[]{\includegraphics[width=0.9\columnwidth]{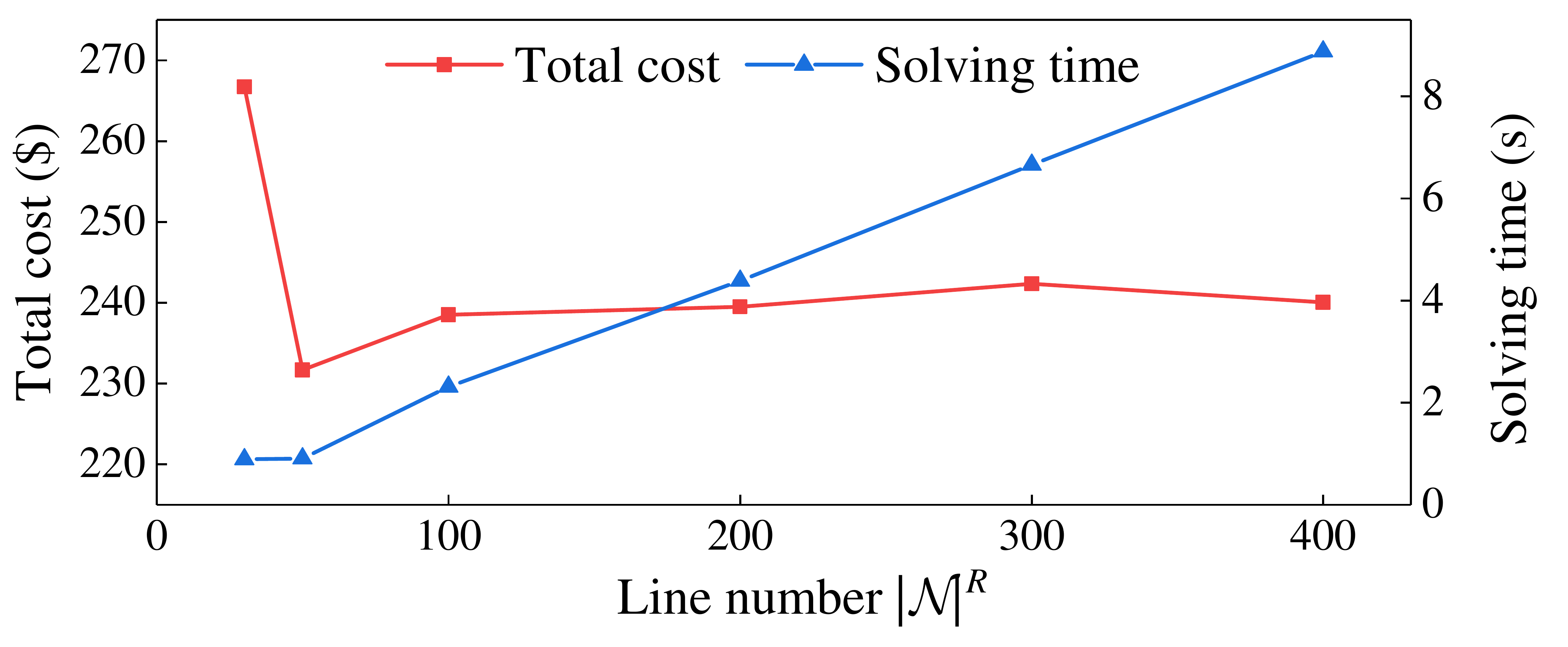}}
	\subfigure[]{\includegraphics[width=0.9\columnwidth]{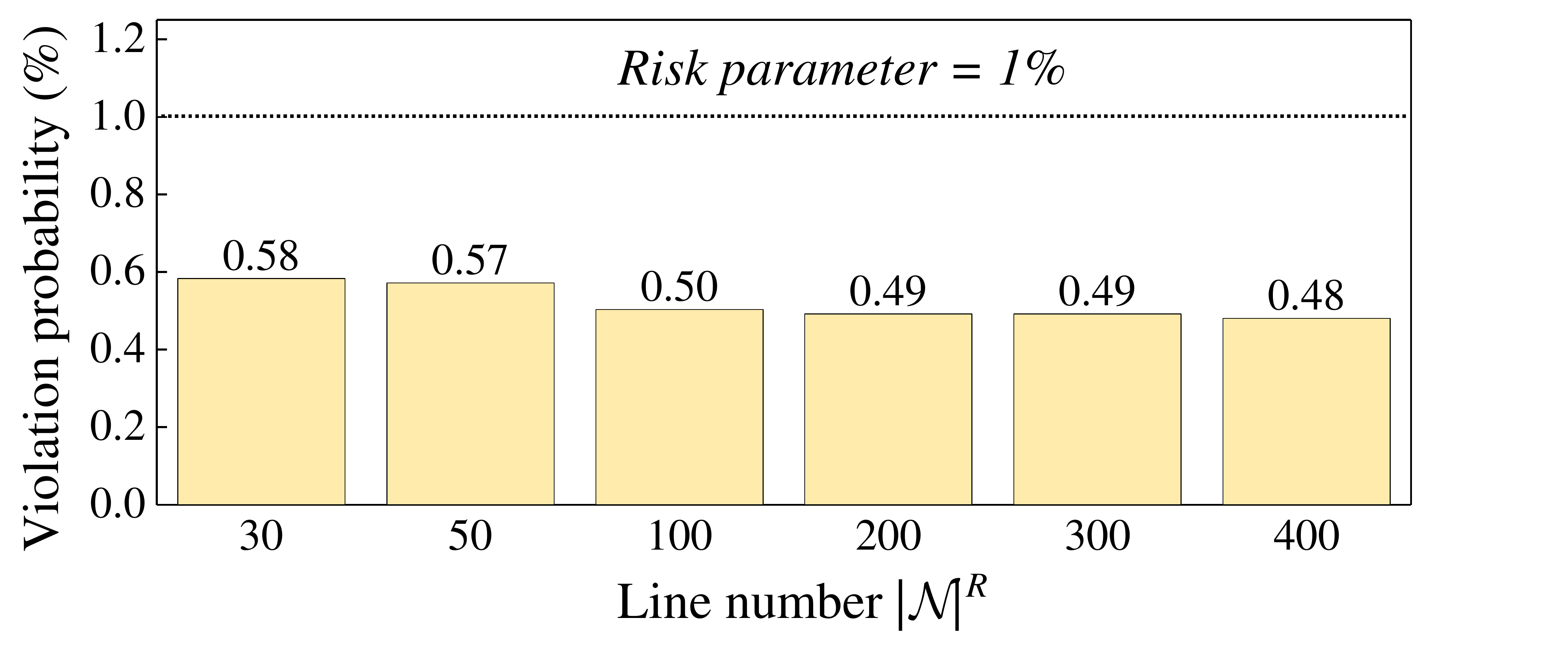}}
	\vspace{1mm}
	\caption{Results of (a) approximation errors caused by the piecewise linearization in  (\ref{eqn_PWL_R}), (b) whole-day total cost and solving time, and (c) maximum violation probability under different line numbers, i.e., $|\mathcal{N}|^R$ with $\epsilon=0.01$. Note the error is defined as $|e^{\rho_2} - R_t^\text{ha}|$.}
	\label{fig_lineNumber}
	\vspace{-4mm}
\end{figure}

In (\ref{eqn_PWL_R}), we employ piecewise linearization to approximate the regulation capacity offer $R_t^\text{ha}$. Fig. \ref{fig_lineNumber} demonstrates the maximum and average approximation errors, whole-day total cost, solving time, and maximum violation probability under different line numbers $|\mathcal{N}|^R$. With the increase of $|\mathcal{N}|^R$, the piecewise linearization used in (\ref{eqn_PWL_R}) becomes more accurate, so both the maximum and average approximation errors decrease, as shown in Fig. \ref{fig_lineNumber}(a). Once $|\mathcal{N}|^R$ reaches 100, the impacts of the approximation error on the optimal solution becomes insignificant. Thus, even if we further increase $|\mathcal{N}|^R$ from 100 to 500, the obtained total cost and maximum violation probability keep almost unchanged, as illustrated in Figs. \ref{fig_lineNumber}(b) and (c). Increasing $|\mathcal{N}|^R$ introduces more binary variables according to (\ref{eqn_PWL_R2}), so the solving time grows rapidly. Nevertheless, in all these cases, the total cost is much less compared to \textbf{B1}, \textbf{B2}, and \textbf{B3}. Moreover, the maximum violation probability is also always lower than the risk parameter. These results further confirm the benefits of the proposed method.

\section{Conclusions} \label{sec_conclusion}
This paper proposes a tractable chance-constrained model to optimize the hour-ahead regulation capacity offers for HVAC systems. It first proposes a temporal compression method to compress the numerous thermodynamic constraints introduced by frequently regulated HVAC power into only a few constraints. Then, a novel mixture-model-based convexification approach is developed to overcome the intractability caused by the non-Gaussian uncertainties from regulation signals. By applying this approach, the chance constraints with these non-Gaussian uncertainties on the LHS can be approximated by tractable SOCP forms with marginal optimality loss.
Numerical experiments confirm that the proposed method can achieve better feasibility performance compared to the widely used Gaussian-assumption-based methods, while its solution is also less conservative than the DRCC methods.

\appendices
\setcounter{table}{0}   

\section{} \label{app_0}
\emph{Proof of Proposition \ref{prop_0}}:
We first define two new vectors $\bm v_1$ and $\bm v_2$, as follows:
\begin{align}
    \bm v_1 = \left[\max_{l \in \mathcal{L}_\tau} f(l), \forall \tau \in \mathcal{T} \right]^{\intercal},
    \bm v_2 = \left[R \min_{l \in \mathcal{L_\tau}}[\bm A \bm s]_{l}, \forall \tau \in \mathcal{T} \right]^{\intercal}. \label{eqn_difine_v}
\end{align}
Then, based on the Minkowski's inequality, we must have:
\begin{align}
     \underbrace{\Vert \bm v_1 + \bm v_2 \Vert_\infty}_{\text{RHS term of  (\ref{eqn_approximation_2})}} \leq \underbrace{\Vert \bm v_1\Vert_\infty + \Vert \bm v_2\Vert_\infty}_{\text{RHS term of  (\ref{eqn_approximation_0})}}. \label{eqn_Minkowski}
\end{align}
By substituting (\ref{eqn_difine_v}) into (\ref{eqn_Minkowski}), we prove the maximum indoor temperature estimated by (\ref{eqn_approximation_2}) is no more than that of (\ref{eqn_approximation_0}). Based on the same way, we can also prove that the minimum indoor temperature estimated by (\ref{eqn_approximation_2}) is no less than that of (\ref{eqn_approximation_0}). This completes the proof.

\section{} \label{app_1}
For the first constraint in (\ref{eqn_thermal_r}), the detail expressions of $\bm \alpha(\bm x)$ and $\beta(\bm x)$ are as follows:
\begin{align}
\text{1st}
\begin{cases}
&\bm \alpha(\bm x) = [a^\text{in}_{\tau}, R^\text{ha}]^\intercal, \quad \bm \omega = [\theta^\text{in}_0, \overline u_\tau]^\intercal,\\
&\beta(\bm x) = \theta^\text{max} - a^\text{out}_{\tau} \theta^\text{out} - a^\text{h}_{\tau} \theta^\text{h} - a^\text{q}_{\tau} p^\text{ha},
\end{cases} \label{eqn_coeff_1}
\end{align}
where 1st represents the first constraint; parameters $a^\text{in}_{\tau}=(a^\text{in})^{\tau\frac{|\mathcal{L}|}{|\mathcal{T}|}}$, $a^\text{out}_{\tau}=a^\text{out} \frac{1-a^\text{in}_{\tau}}{1-a^\text{in}}$, $a^\text{h}_{\tau}=a^\text{h} \frac{1-a^\text{in}_{\tau}}{1-a^\text{in}}$, and $a^\text{q}_{\tau}=a^\text{q} \frac{1-a^\text{in}_{\tau}}{1-a^\text{in}}$. Similarly, for the rest chance constraints, we have
\begin{align}
&\text{2nd}
\begin{cases}
&\bm \alpha(\bm x) = [a^\text{in}_{\tau+1}, R^\text{ha}]^\intercal, \quad \bm \omega = [\theta^\text{in}_0, \overline u_\tau]^\intercal,\\
&\beta(\bm x) = \theta^\text{max} - a^\text{out}_{\tau+1} \theta^\text{out} - a^\text{h}_{\tau+1} \theta^\text{h} - a^\text{q}_{\tau+1} p^\text{ha},
\end{cases}\\
&\text{3rd}
\begin{cases}
&\bm \alpha(\bm x) = [a^\text{in}_{\tau}, R^\text{ha}]^\intercal, \quad \bm \omega = [-\theta^\text{in}_0, -\underline u_\tau]^\intercal,\\
&\beta(\bm x) = a^\text{out}_{\tau} \theta^\text{out} + a^\text{h}_{\tau} \theta^\text{h} + a^\text{q}_{\tau} p^\text{ha} - \theta^\text{min},
\end{cases}\\
&\text{4th}
\begin{cases}
&\bm \alpha(\bm x) = [a^\text{in}_{\tau+1}, R^\text{ha}]^\intercal, \quad \bm \omega = [-\theta^\text{in}_0, -\underline u_\tau]^\intercal,\\
&\beta(\bm x) = a^\text{out}_{\tau+1} \theta^\text{out} + a^\text{h}_{\tau+1} \theta^\text{h} + a^\text{q}_{\tau+1} p^\text{ha} - \theta^\text{min}.
\end{cases} \label{eqn_coeff_4}
\end{align}
Note that only one single variable, i.e., $R^\text{ha}$, is contained in the four $\bm \alpha(\bm x)$ in (\ref{eqn_coeff_1})-(\ref{eqn_coeff_4}).

\section{} \label{app_2}
\emph{Proof of Proposition \ref{prop_1}}:
The convexity of function $\ln{\Phi^{-1}(y_j)}$ can be analyzed by its second-order derivative:
\begin{align}
\frac{d^2\ln{\Phi^{-1}(y_j)}}{(dy_j)^2} = \frac{-\phi(v_j)-v_j \phi'(v_j)}{v_j^2(\phi(v_j))^3}, \quad \forall j \in \mathcal{J}, \label{eqn_derivative}
\end{align}
where $v_j=\Phi^{-1}(y_j)$; $\phi(\cdot)$ is the PDF of the standard normal distribution and $\phi'(\cdot)$ is its first-order derivative. Note the PDF $\phi(\cdot)$ is always nonnegative. Considering that a convex function has a nonnegative second-order derivative, we can get the convex condition for function $\ln{\Phi^{-1}(y_j)}$, as follows:
\begin{align}
&\frac{-\phi(v_j)-v_j \phi'(v_j)}{v_j^2(\phi(v_j))^3} \geq 0 \Leftrightarrow \phi(v_j)+v_j \phi'(v_j) \leq 0 \notag\\
&\Leftrightarrow \frac{e^{-v_j^2/2}}{\sqrt{2\pi}} - v_j^2 \frac{e^{-v_j^2/2}}{\sqrt{2\pi}} \leq 0  \Leftrightarrow v_j^2 \geq 1.
\end{align}
According to the definition of chance constraints, we have $y_j\geq 0.5$. Thus, variable $v_j$ is non-negative. As a result, the above inequality can be further converted into:
\begin{align}
v_j^2 \geq 1 \Leftrightarrow v_j \geq 1 \Leftrightarrow y_j \geq \Phi(1).
\end{align}
Similarly, its concave condition can be obtained by:
\begin{align}
&\frac{-\phi(v_j)-v_j \phi'(v_j)}{v_j^2(\phi(v_j))^3} \leq 0 \Leftrightarrow v_j^2 \leq 1 \Leftrightarrow y_j \leq \Phi(1).
\end{align}

We uniformly select $N+1$ points on function $\ln{\Phi^{-1}(y_j)}$ (recorded by $y_j^{(n)}, \forall n \in \mathcal{N}$) and let $y_j^{(0)}=\Phi(1)<y_j^{(1)}<\cdots<y_j^{(N)}$. With points $y_j^{(1)},\cdots,y_j^{(N)}$, N-1 line segments can be constructed by connecting these points in sequence (denoted by $\lambda_n y_j + \gamma_n, \forall n \in \mathcal{N}/\{0\}$). Since function $\ln{\Phi^{-1}(y_j)}$ is convex when $y_j\geq \Phi(1)$, according to the definition of convex functions, we must have:
\begin{align}
\ln{\Phi^{-1}(y_j)} \leq \max_{n \in \mathcal{N}/\{0\}}{\{\lambda_n y_j + \gamma_n}\}, \text{for } y_j \geq \Phi(1). \label{eqn_proof_1}
\end{align}
For the region $y_j \leq \Phi(1)$, function $\ln{\Phi^{-1}(y_j)}$ is concave. According to the concave function's first-order condition, function $\ln{\Phi^{-1}(y_j)}$ should always be equal to or below its tangent. Thus, by letting $\lambda_0 y_j + \gamma_0$ as the tangent of function $\ln{\Phi^{-1}(y_j)}$ at $y_j = \Phi(1)$, we have:
\begin{align}
\ln{\Phi^{-1}(y_j)} \leq \lambda_0 y_j + \gamma_0, \quad \text{for } y_j \leq \Phi(1). \label{eqn_proof_2}
\end{align}
By combining (\ref{eqn_proof_1}) and (\ref{eqn_proof_2}), we prove \textbf{Proposition} \ref{prop_1}.

\footnotesize
\bibliographystyle{ieeetr}
\bibliography{ref}
\end{document}